# THE LARGEST SAMPLE EIGENVALUE DISTRIBUTION IN THE RANK 1 QUATERNIONIC SPIKED MODEL OF WISHART ENSEMBLE

By Dong Wang

*Brandeis University*

We solve the largest sample eigenvalue distribution problem in the rank 1 spiked model of the quaternionic Wishart ensemble, which is the first case of a statistical generalization of the Laguerre symplectic ensemble (LSE) on the soft edge. We observe a phase change phenomenon similar to that in the complex case, and prove that the new distribution at the phase change point is the GOE Tracy–Widom distribution.

**1. Introduction.** The Wishart ensemble is defined as follows [24]:

Consider $M$ independent observation $\mathbf{x}_1, \ldots, \mathbf{x}_M$ of an $N$-variate normal distribution with mean 0 and covariance matrix $\Sigma$. Here the values of the normal distribution can be real, complex or even quaternion. If the variables are complex or quaternionic, then the definition of the mean is as usual, and the (co)variance is defined as

$$\mathrm{cov}(x,y) = \mathrm{E}((x-\bar{x})(y-\bar{y})^*),$$

where $\bar{x}$ (resp. $\bar{y}$) is the mean of $x$ (resp. $y$), and $*$ is the complex or quaternionic conjugation operator. Then $\Sigma$ is a real symmetric/Hermitian/quaternionic Hermitian matrix. Without loss of generality, we assume $\Sigma$ to be a diagonal matrix, with *population eigenvalues* $l = (l_1, \ldots, l_N)$. If we put the above data into an $N \times M$ double array $\mathbf{X} = (\mathbf{x}_1 : \cdots : \mathbf{x}_M)$, then the positively defined real symmetric/Hermitian/quaternionic hermitian matrix $S = \frac{1}{M} \mathbf{X} \mathbf{X}^*$ is the *sample matrix* and its eigenvalues $\lambda = (\lambda_1, \ldots, \lambda_N)$ are *sample eigenvalues*. ($\mathbf{X}^*$ is the transpose, Hermitian transpose or quaternionic Hermitian transpose of $\mathbf{X}$ depending on type of $\mathbf{X}$'s entries.) The probability space of $\lambda_i$'s is called the Wishart ensemble.

It is a classical result [2] that (in the real category) if $M \gg N$, $\lambda_i$'s are good approximations of $l_i$'s. But if $M$ and $N$ are of the same order of magnitude, that is, $M/N = \gamma^2 \geq 1$ and $M$ and $N$ are very large, the problem is subtler. The simplest case with $\Sigma = I$, the *white* Wishart ensemble, is the Laguerre ensemble, well studied in random matrix theory (RMT) under the name LOE, LUE and LSE—they are abbreviations of Laguerre Orthogonal/Unitary/Symplectic Ensemble, and GOE, GUE and GSE appearing later are abbreviations of Gaussian Orthogonal/Unitary/Symplectic Ensemble—over all the three base fields, respectively.

Naturally, the next question is: If $\Sigma$ is slightly deviate from $I$, such that $l_i = 1 + a_i$, $i = 1, \ldots, k$, and $l_{k+1} = \cdots = l_N = 1$, what is the distribution of the $\lambda_i$'s? This is called the *spiked model* [19] and $k$ is defined as its rank.

If $M$ and $N$ are very large and $k$ and $a_i$'s are small constants, the density of $\lambda_i$'s is the same as that in the white Wishart model, proved in [22] in real and complex categories. The distribution of the largest sample eigenvalue, however, may change. For the complex ensemble, Baik, Ben Arous and Péché [4] solved the problem completely. They show that if $\max(a_i)$ is smaller than a threshold, then the distribution of the largest sample eigenvalue is the same as that in the white ensemble, which is the GUE Tracy–Widom distribution, but if $\max(a_i)$ exceeds the threshold, that distribution is changed into a Gaussian whose mean and variance depend on $\max(a_i)$. Furthermore, in the case that $\max(a_i)$ equals the critical value, they find a series of new distributions, indexed by the multiplicity of $\max(a_i)$.

In the real category, which is practically the most important and mathematically the most difficult, much less is known. In this paper I solve the distribution of the largest sample eigenvalue for the rank 1 spiked model in the quaternionic category. I believe the similarity of LOE and LSE [13] suggests that the solution to the quaternionic spiked model is an intermediate step toward the solution to the real one.

1.1. *Some known results for the largest sample eigenvalue in white and rank 1 spiked models.* In latter part of the paper, we concentrate on the distribution of the largest sample eigenvalue in the rank 1 spiked model, so denote $a$ to be the only perturbation parameter.

The result in the complex category is complete. First we recall the result for the complex white Wishart ensemble.

PROPOSITION 1. *The distribution of the largest sample eigenvalue in the complex white Wishart ensemble satisfies that,* $\max(\lambda)$ *almost surely approaches [15]* $(1+\gamma^{-1})^2$ *with fluctuation scale $M^{-2/3}$, and [11, 18]*

$$\lim_{M \to \infty} \mathbb{P}\bigg((\max(\lambda) - (1+\gamma^{-1})^2) \cdot \frac{\gamma M^{2/3}}{(1+\gamma)^{4/3}} \leq T\bigg) = F_{\mathrm{GUE}}(T),$$

*where $F_{\mathrm{GUE}}$ is the GUE Tracy–Widom distribution.*

The GUE Tracy–Widom distribution is defined by Fredholm determinant [11, 29]:

$$F_{\text{GUE}}(T) = \det(1 - K_{\text{Airy}}(\xi,\eta)\chi_{(T,\infty)}(\eta)),$$

where $\chi_{(T,\infty)}$ is the step function:

$$\chi_{(T,\infty)}(\eta) = \begin{cases} 1, & \text{if } \eta \in (T,\infty), \\ 0, & \text{otherwise,} \end{cases}$$

and $K_{\text{Airy}}(\xi,\eta)$ is the well-known Airy kernel defined by the Airy function $\text{Ai}(x)$:

(1) $$K_{\text{Airy}}(\xi,\eta) = \int_0^\infty \text{Ai}(\xi+t)\,\text{Ai}(\eta+t)\,dt.$$

The Airy function can be defined in different ways, and here we take an integral representation suitable for our asymptotic analysis [4]:

(2) $$\text{Ai}(\xi) = \frac{-1}{2\pi i}\int_{\Gamma^\infty} e^{-\xi z + 1/3 z^3}\,dz,$$

where $\Gamma^\infty = \Gamma_1^\infty \cup \Gamma_2^\infty \cup \Gamma_3^\infty$, which are defined as (see Figure 1)

$$\Gamma_1^\infty = \{-te^{\pi i/3}|-\infty < t \leq -1\}, \qquad \Gamma_2^\infty = \{e^{-t\pi i}|-\tfrac{1}{3} \leq t \leq \tfrac{1}{3}\},$$

$$\Gamma_3^\infty = \{te^{5\pi i/3}|1 \leq t < \infty\}.$$

The breakthrough in the complex category is by [4], which is for any finite rank spiked model. In the rank 1 case, it is:

PROPOSITION 2. *In the rank 1 complex spiked model:*

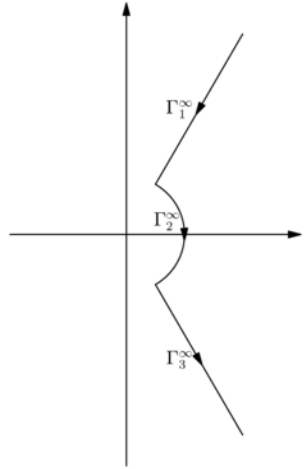

FIG. 1. $\Gamma^\infty$.

1. If $-1 < a < \gamma^{-1}$, then the distribution of the largest sample eigenvalue is the same as that of the complex white Wishart ensemble in Proposition 1.
2. If $a = \gamma^{-1}$, then the limit and the fluctuation scale are the same as those of the complex white Wishart ensemble, but the distribution function is

$$\lim_{M \to \infty} \mathbb{P}\left( (\max(\lambda) - (1+\gamma^{-1})^2) \cdot \frac{\gamma M^{2/3}}{(1+\gamma)^{4/3}} \leq T \right) = F_{\text{GUE1}}(T). \tag{3}$$

3. If $a > \gamma^{-1}$, then the limit and the fluctuation scale are changed as well as the distribution function, which is a Gaussian:

$$\lim_{M \to \infty} \mathbb{P}\left( \left( \max(\lambda) - (a+1)\left(1 + \frac{1}{\gamma^2 a}\right) \right) \cdot \frac{\sqrt{M}}{(a+1)\sqrt{1 - 1/(\gamma^2 a^2)}} \leq T \right)$$
$$= \int_{-\infty}^{T} \frac{1}{\sqrt{2\pi}} e^{-t^2/2} \, dt. \tag{4}$$

The function $F_{\text{GUE1}}$ occurring in (3) is defined similarly to $F_{\text{GUE}}$ [4]:

$$F_{\text{GUE1}}(T) = \det(1 - (K_{\text{Airy}}(\xi, \eta) + s^{(1)}(\xi) \operatorname{Ai}(\eta)) \chi_{(T,\infty)}(\eta)), \tag{5}$$

where $s^{(1)}$ is one of a series of functions defined in [4], and has the integral representation

$$s^{(1)}(\eta) = \frac{1}{2\pi i} \int_{\bar{\Gamma}^\infty} e^{-\eta z + 1/3 z^3} \frac{1}{z} \, dz \quad \text{and} \quad s^{(1)}(\eta) = 1 - \int_\eta^\infty \operatorname{Ai}(t) \, dt,$$

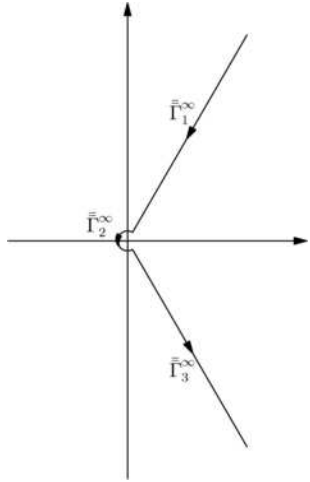

Fig. 2. $\bar{\bar{\Gamma}}^\infty$.

where $\bar{\bar{\Gamma}}^\infty = \bar{\bar{\Gamma}}_1^\infty \cup \bar{\bar{\Gamma}}_2^\infty \cup \bar{\bar{\Gamma}}_3^\infty$, which are defined as (see Figure 2 $\varepsilon$ is a positive constant, used later)

$$\bar{\bar{\Gamma}}_1^\infty = \left\{-te^{\pi i/3}\Big|-\infty < t \leq -\frac{\varepsilon}{2}\right\}, \qquad \bar{\bar{\Gamma}}_2^\infty = \left\{\frac{\varepsilon}{2}e^{t\pi i}\Big|\frac{1}{3} \leq t \leq \frac{5}{3}\right\},$$

$$\bar{\bar{\Gamma}}_3^\infty = \left\{te^{5\pi i/3}\Big|\frac{\varepsilon}{2} \leq t < \infty\right\}.$$

REMARK 1. The kernel in (5) is not in trace class, but the Fredholm determinant is well defined and we can easily conjugate it into a trace class kernel. Several kernels below are in similar situations.

In the real category, we have the result for the real white Whishart ensemble:

PROPOSITION 3. *The distribution of the largest sample eigenvalue in the real white Wishart ensemble satisfies that, $\max(\lambda)$ almost surely approaches* [15] $(1+\gamma^{-1})^2$ *with fluctuation scale $M^{-2/3}$, and* [19]

$$\lim_{M\to\infty} \mathbb{P}\bigg((\max(\lambda) - (1+\gamma^{-1})^2) \cdot \frac{\gamma M^{2/3}}{(1+\gamma)^{4/3}} \leq T\bigg) = F_{\text{GOE}}(T),$$

*where $F_{\text{GOE}}$ is the GOE Tracy–Widom distribution.*

Here the function $F_{\text{GOE}}$ is defined by the Fredholm determinant of a matrix integral operator [30]:

$$F_{\text{GOE}}(T) = \sqrt{\det(I - P_{\text{GOE}}(\xi, \eta))}$$

and

$$P_{\text{GOE}}(\xi, \eta) = \chi_{(T,\infty)}(\xi) \begin{pmatrix} S_1(\xi, \eta) & SD_1(\xi, \eta) \\ IS_1(\xi, \eta) - \frac{1}{2}\operatorname{sgn}(x-y) & S_1(\eta, \xi,) \end{pmatrix} \chi_{(T,\infty)}(\eta),$$

where

$$S_1(\xi, \eta) = K_{\text{Airy}}(\xi, \eta) - \frac{1}{2}\operatorname{Ai}(\xi)\int_\eta^\infty \operatorname{Ai}(t)\,dt + \frac{1}{2}\operatorname{Ai}(\xi),$$

$$SD_1(\xi, \eta) = -\frac{\partial}{\partial \eta}K_{\text{Airy}}(\xi, \eta) - \frac{1}{2}\operatorname{Ai}(\xi)\operatorname{Ai}(\eta),$$

$$IS_1(\xi, \eta) = -\int_\xi^\infty K_{\text{Airy}}(t, \eta)\,dt + \frac{1}{2}\int_\xi^\infty \operatorname{Ai}(t)\,dt \int_\eta^\infty \operatorname{Ai}(t)\,dt$$

$$- \frac{1}{2}\int_\xi^\infty \operatorname{Ai}(t)\,dt + \frac{1}{2}\int_\eta^\infty \operatorname{Ai}(t)\,dt.$$

REMARK 2. We have a more convenient form of $F_{\text{GOE}}$ [12]:

$$F_{\text{GOE}} = \sqrt{\det(1 - (K_{\text{Airy}}(\xi,\eta) + s^{(1)}(\xi)\operatorname{Ai}(\eta))\chi_{(T,\infty)}(\eta))}, \tag{6}$$

so [4]

$$F_{\text{GUE1}}(T) = (F_{\text{GOE}}(T))^2.$$

In the real spiked model, Baik and Silverstein [8] compute the almost sure limit of the largest population eigenvalue, which is the same as that in the complex category, and Paul [25] proves the Gaussian distribution property in the case $a > \gamma^{-1}$, which is similar to (4). Neither of their methods can find the distribution function when $a \leq \gamma^{-1}$.

For the quaternionic white Wishart ensemble, we have:

PROPOSITION 4. *The distribution of the largest sample eigenvalue in the quaternionic white Wishart ensemble satisfies that,* $\max(\lambda)$ *almost surely approaches* $(1+\gamma^{-1})^2$ *with fluctuation scale* $M^{-2/3}$, *and* [14]

$$\lim_{M\to\infty} \mathbb{P}\left((\max(\lambda) - (1+\gamma^{-1})^2)\cdot \frac{\gamma(2M)^{2/3}}{(1+\gamma)^{4/3}} \leq T\right) = F_{\text{GSE}}(T),$$

*where* $F_{\text{GSE}}$ *is the GSE Tracy–Widom distribution.*

Here the function $F_{\text{GSE}}$ is defined by the Fredholm determinant of a matrix integral operator [30]:

$$F_{\text{GSE}}(T) = \sqrt{\det(I - \widehat{P}(\xi,\eta))}$$

and

$$\hat{P}(\xi,\eta) = \chi_{(T,\infty)}(\xi)\begin{pmatrix} \widehat{S}_4(\xi,\eta) & \widehat{SD}_4(\xi,\eta) \\ \widehat{IS}_4(\xi,\eta) & \widehat{S}_4(\eta,\xi,) \end{pmatrix}\chi_{(T,\infty)}(\eta),$$

where

$$\widehat{S}_4(\xi,\eta) = \frac{1}{2}K_{\text{Airy}}(\xi,\eta) - \frac{1}{4}\operatorname{Ai}(\xi)\int_\eta^\infty \operatorname{Ai}(t)\,dt,$$

$$\widehat{SD}_4(\xi,\eta) = -\frac{1}{2}\frac{\partial}{\partial \eta}K_{\text{Airy}}(\xi,\eta) - \frac{1}{4}\operatorname{Ai}(\xi)\operatorname{Ai}(\eta),$$

$$\widehat{IS}_4(\xi,\eta) = -\frac{1}{2}\int_\xi^\infty K_{\text{Airy}}(t,\eta)\,dt + \frac{1}{4}\int_\xi^\infty \operatorname{Ai}(t)\,dt\int_\eta^\infty \operatorname{Ai}(t)\,dt.$$

1.2. *Statement of main results.* The main theorem in this paper is:

THEOREM 1. *In the rank 1 quaternionic spiked model:*

1. *If $-1 < a < \gamma^{-1}$, then the distribution of the largest sample eigenvalue is the same as that of the quaternionic white Wishart ensemble in Proposition 4.*
2. *If $a = \gamma^{-1}$, then the limit and the fluctuation scale are the same as those of the quaternionic white Wishart ensemble, but the distribution function is*

$$\lim_{M \to \infty} \mathbb{P}\left(\left(\max(\lambda) - \left(\frac{\gamma+1}{\gamma}\right)^2\right) \cdot \frac{\gamma(2M)^{2/3}}{(1+\gamma)^{4/3}} \leq T\right) = F_{\text{GSE1}}(T).$$

3. *If $a > \gamma^{-1}$, then the limit and the fluctuation scale are changed as well as the distribution function, which is a Gaussian:*

$$\lim_{M \to \infty} \mathbb{P}\left(\left(\max(\lambda) - (a+1)\left(1 + \frac{1}{\gamma^2 a}\right)\right) \cdot \frac{\sqrt{2M}}{(a+1)\sqrt{1 - 1/(\gamma^2 a^2)}} \leq T\right)$$
$$= \int_{-\infty}^T \frac{1}{\sqrt{2\pi}} e^{-t^2/2} \, dt.$$

Here the function $F_{\text{GSE1}}$ is defined by the Fredholm determinant of a matrix integral operator:

$$F_{\text{GSE1}}(T) = \sqrt{\det(I - \overline{\overline{P}}(\xi, \eta))}$$

and

$$\overline{\overline{P}}(\xi, \eta) = \chi_{(T,\infty)}(\xi) \begin{pmatrix} \overline{\overline{S}}_4(\xi, \eta) & \overline{\overline{SD}}_4(\xi, \eta) \\ \overline{\overline{IS}}_4(\xi, \eta) & \overline{\overline{S}}_4(\eta, \xi, ) \end{pmatrix} \chi_{(T,\infty)}(\eta),$$

where

$$\overline{\overline{S}}_4(\xi, \eta) = \widehat{S}_4(\xi, \eta) + \frac{1}{2} \text{Ai}(\xi), \qquad \overline{\overline{SD}}_4(\xi, \eta) = \widehat{SD}_4(\xi, \eta),$$

$$\overline{\overline{IS}}_4(\xi, \eta) = \widehat{IS}_4(\xi, \eta) - \frac{1}{2} \int_\xi^\infty \text{Ai}(t) \, dt + \frac{1}{2} \int_\eta^\infty \text{Ai}(t) \, dt.$$

Although the distribution $F_{\text{GSE1}}$ seems to be new, we have that

THEOREM 2.

$$F_{\text{GSE1}}(T) = F_{\text{GOE}}(T).$$

1.3. *Relation with other models and conjecture on the rank 1 real spiked model.* The results of Theorems 1 and 2 give a phase transition pattern $F_{\text{GSE}}$–$F_{\text{GOE}}$–Gaussian as the parameter $a$ increases from $-1$ to $+\infty$. This pattern appears as limiting distributions indexed by a parameter in several other combinatorial and statistical physical models, for example, the lengths of the longest monotone subsequences of random involutions with condition on the number of fixed points [6] and the symmetrized last passage percolation [7] studied by Baik and Rains. In semi-infinite totally asymmetric simple exclusion process [26] studied by Prähofer and Spohn, and the symmetric polynuclear growth process [5] studied by Baik et al., 2-dimensional phase transition diagrams are obtained, and the 1-dimensional $F_{\text{GSE}}$–$F_{\text{GOE}}$–Gaussian pattern is contained in both of them.

Although there is no model which can give hints to the rank 1 real spiked model, it is plausible that it has a phase transition from $F_{\text{GOE}}$ to Gaussian for the limiting distributions of the largest sample eigenvalue as $a$ goes across $\gamma^{-1}$. Based on the duality of orthogonal and symplectic models from the Virasoro structure's point of view, we have:

CONJECTURE 1. In the rank 1 real spiked model:

1. If $-1 < a < \gamma^{-1}$, then the distribution of the largest sample eigenvalue is the same as that of the real white Wishart ensemble in Proposition 3.
2. If $a = \gamma^{-1}$, then the limit and the fluctuation scale are the same as those of the quaternionic white Wishart ensemble, but the distribution function is

$$\lim_{M \to \infty} \mathbb{P}\left(\left(\max(\lambda) - \left(\frac{\gamma+1}{\gamma}\right)^2\right) \cdot \frac{\gamma(M/2)^{2/3}}{(1+\gamma)^{4/3}} \leq T\right) = F_{\text{GSE}}(T).$$

3. If $a > \gamma^{-1}$, then the limit and the fluctuation scale are changed as well as the distribution function, which is a Gaussian (proved by Paul in [25]):

$$\lim_{M \to \infty} \mathbb{P}\left(\left(\max(\lambda) - (a+1)\left(1 + \frac{1}{\gamma^2 a}\right)\right) \cdot \frac{\sqrt{M/2}}{(a+1)\sqrt{1 - 1/(\gamma^2 a^2)}} \leq T\right)$$
$$= \int_{-\infty}^{T} \frac{1}{\sqrt{2\pi}} e^{-t^2/2}\, dt.$$

1.4. *Structure of the paper.* In Section 2 we use combinatorial techniques to express the joint distribution function of $\{\lambda_j\}$, and then by skew orthogonal polynomial techniques express the distribution function of $\max(\lambda_j)$ in the square root of a Fredholm determinant of a matrix integral operator. In Section 3 we do asymptotic analysis on the kernel of the matrix integral operator, and prove the three cases of Theorem 1 in the three subsections, respectively. Section 4 contains the proof of Theorem 2. In the proof of Theorem 1, we use some trace norm convergence results which generalize the old result on the LUE [11], and we give a method of proof to them in the Appendix.

## 2. The Fredholm determinantal formula.

2.1. *The joint distribution function.* In this subsection, we prove the following:

THEOREM 3. *The joint probability distribution function of $\lambda$ in the quaternionic spiked model is*

$$(7) \qquad P(\lambda) = \frac{1}{C} \tilde{V}^4(\lambda) \prod_{j=1}^{N} (\lambda_j^{2(M-N)+1} e^{-2M\lambda_j}).$$

In this paper, $C$ stands for any constants, and here

$$\tilde{V}^4(\lambda) = \begin{vmatrix} 1 & 0 & \cdots & 1 & 0 \\ \lambda_1 & 1 & \cdots & \lambda_N & 1 \\ \lambda_1^2 & 2\lambda_1 & \cdots & \lambda_N^2 & 2\lambda_N \\ \vdots & \vdots & \cdots & \vdots & \vdots \\ \lambda_1^{2N-2} & (2N-2)\lambda_1^{2N-3} & \cdots & \lambda_N^{2N-2} & (2N-2)\lambda_N^{2N-3} \\ e^{a/(1+a)2M\lambda_1} & \frac{a}{1+a}2Me^{a/(1+a)2M\lambda_1} & \cdots & e^{a/(1+a)2M\lambda_N} & \frac{a}{1+a}2Me^{a/(1+a)2M\lambda_N} \end{vmatrix},$$

the determinant of a $2N \times 2N$ matrix whose $(2N, 2k-1)$ entry is $e^{a/(1+a)2M\lambda_k}$, $(j, 2k-1)$ entry is $\lambda_k^{j-1}$ for $j = 1, \ldots, 2N-1$, and $2i$th column is the derivative of the $(2i-1)$st column. $\tilde{V}^4(\lambda)$ is a variation of the $V(\lambda)^4$ appearing in the LSE (see [23] and (9)).

For the Wishart ensemble defined in the introduction section, we first have the distribution function for the sample matrix in the $N \times N$ positive definite quaternionic Hermitian matrix space [3]:

$$P(S) = \frac{1}{C} e^{-2M\Re \operatorname{Tr}(\Sigma^{-1}S)} (\det S)^{2(M-N)+1}.$$

REMARK 3. Due to the noncommutativity of the quaternions, $\det S$ is not well defined in the usual way. Since $S$ is quaternionic Hermitian, we can diagonalize it into a real-valued diagonal matrix by the conjugation of a quaternionic unitary matrix $U$, and define

$$\det S = \prod_N \text{eigenvalues of } USU^*.$$

REMARK 4. In the distribution function in real and complex categories of sample matrices, we do not need to take the real part of the trace, since the trace is already real. Unfortunately, this does not hold in the quaternionic category due to its noncommutativity, and luckily $\Re \operatorname{Tr}$ behaves better. [For example, $\Re \operatorname{Tr}(AB) = \Re \operatorname{Tr}(BA)$, but $\operatorname{Tr}(AB) \neq \operatorname{Tr}(BA)$ in general.]

The distribution function for sample eigenvalues $\lambda$, the eigenvalues of $S$, is

$$(8) \quad P(\lambda) = \frac{1}{C}(V(\lambda))^4 \prod_{j=1}^{N} \lambda_j^{2(M-N)+1} \int_{Q \in Sp(N)} e^{-2M\Re \operatorname{Tr}(\Sigma^{-1}Q\Lambda Q^{-1})} \, dQ,$$

where we integrate on the compact symplectic group with the Haar measure, $V(\lambda) = \prod_{i<j}(\lambda_i - \lambda_j)$ is the Vandermonde, and $\Lambda = \operatorname{diag}(\lambda_1, \ldots, \lambda_N)$. (See [23] for a derivation of the similar GSE case.)

If the perturbation parameter $a = 0$, then $l_1 = l_2 = \cdots = l_N = 1$,

$$\int_{Q \in Sp(N)} e^{-2M\Re \operatorname{Tr}(\Sigma^{-1}Q\Lambda Q^{-1})} \, dQ = \prod_{j=1}^{N} e^{-2M\lambda_j}$$

and

$$(9) \qquad P(\lambda) = \frac{1}{C}(V(\lambda))^4 \prod_{j=1}^{N} (\lambda_j^{2(M-N)+1} e^{-2M\lambda_j}),$$

is the standard LSE [23].

Generally,

$$(10) \quad \begin{aligned} &\int_{Q \in Sp(N)} e^{-2M\Re \operatorname{Tr}(\Sigma^{-1}Q\Lambda Q^{-1})} \, dQ \\ &= \int_{Q \in Sp(N)} e^{-2M\Re \operatorname{Tr}(IQ\Lambda Q^{-1})} e^{-2M\Re \operatorname{Tr}((\Sigma^{-1}-I)Q\Lambda Q^{-1})} \, dQ \\ &= \prod_{j=1}^{N} e^{-2M\lambda_j} \int_{Q \in Sp(N)} e^{2M\Re \operatorname{Tr}((I-\Sigma^{-1})Q\Lambda Q^{-1})} \, dQ. \end{aligned}$$

Then by the integral formula of the quaternionic Zonal polynomials [17], we get

$$(11) \quad \begin{aligned} &\int_{Q \in Sp(N)} e^{2M\Re \operatorname{Tr}((I-\Sigma^{-1})Q\Lambda Q^{-1})} \, dQ \\ &= \sum_{j=0}^{\infty} \frac{(2M)^j}{j!} \sum_{\substack{l(\kappa) \leq N \\ \kappa \vdash j}} \frac{C_\kappa^{(1/2)}(I-\Sigma^{-1}) C_\kappa^{(1/2)}(\Lambda)}{C_\kappa^{(1/2)}(I_N)}, \end{aligned}$$

where $C_\kappa^{(1/2)}(x_1, \ldots, x_N)$ is the $N$ variable quaternionic Zonal polynomial, that is, the Jack polynomial with the parameter $\alpha = 1/2$ (see [21] and [27]) and the $C$-normalization [10], so that [$\kappa = (k_1, \ldots, k_l)$, $k_1 \geq k_2 \geq \cdots \geq k_l > 0$, then $l(\kappa) = l$]

$$\sum_{\substack{l(\kappa) \leq m \\ \kappa \vdash k}} C_\kappa^{(1/2)}(x_1, \ldots, x_m) = (x_1 + \cdots + x_m)^k.$$

In the formula, a symmetric polynomial of a matrix is equivalent to the symmetric polynomial of its eigenvalues, so

$$C_\kappa^{(1/2)}(I - \Sigma^{-1}) = C_\kappa^{(1/2)}\left(\frac{a}{1+a}, 0, \ldots, 0\right).$$

Since all variables except for one vanish in $C_\kappa^{(1/2)}(I - \Sigma^{-1})$, we simply find

(12) $$C_\kappa^{(1/2)}(I - \Sigma^{-1})|_{l(\kappa)>1} = 0.$$

We have [27]

$$C_{(j)}^{(1/2)}\left(\frac{a}{1+a}, 0, \ldots, 0\right) = \left(\frac{a}{1+a}\right)^j$$

and since the number of variables is $N$ [27]

$$C_{(j)}^{(1/2)}(1, \ldots, 1) = \frac{1}{(j+1)!} \prod_{i=0}^{j-1}(2N+i),$$

so with (11) and (12), we get

$$\int_{Q \in Sp(N)} e^{2M \Re \operatorname{Tr}((I-\Sigma^{-1})Q\Lambda Q^{-1})} \, dQ$$

$$= \sum_{j=0}^{\infty} \frac{(2M)^j}{j!} \frac{C_{(j)}^{(1/2)}(I - \Sigma^{-1}) C_{(j)}^{(1/2)}(\Lambda)}{C_{(j)}^{(1/2)}(I_N)}$$

$$= \sum_{j=0}^{\infty} \frac{j+1}{\prod_{i=0}^{j-1}(2N+i)} \left(\frac{a}{1+a} 2M\right)^j C_{(j)}^{(1/2)}(\Lambda).$$

In [27] there is an identity

$$\sum_{j=0}^{\infty}(j+1)C_{(j)}^{(1/2)}(\Lambda)t^j = \prod_{j=1}^{N} \frac{1}{(1-\lambda_j t)^2}.$$

Comparing it with the well-known identity for Schur polynomials

$$\sum_{j=0}^{\infty} s_{(j)}(\Lambda)t^j = \prod_{j=1}^{N} \frac{1}{1-\lambda_j t},$$

we get the identity

(13) $$(j+1)C^{(1/2)}_{(j)}(\Lambda) = s_{(j)}(\lambda_1, \lambda_1, \lambda_2, \lambda_2, \ldots, \lambda_N, \lambda_N),$$

with each $\lambda_i$ appearing twice as variables of the $s_{(j)}$. For notational simplicity, we denote the right-hand side of (13) as $\tilde{s}_{(j)}(\Lambda)$, which is a plethysm [21]

$$\tilde{s}_{(j)}(\Lambda) = s_{(j)} \circ 2p_1(\Lambda).$$

Now we get

(14) $$\int_{Q \in Sp(N)} e^{2M \Re \operatorname{Tr}((I-\Sigma^{-1})Q\Lambda Q^{-1})} dQ$$
$$= \sum_{j=0}^{\infty} \frac{1}{\prod_{i=0}^{j-1}(2N+i)} \left(\frac{a}{1+a} 2M\right)^j \tilde{s}_{(j)}(\Lambda).$$

Then we need a lemma to simplify (14) further.

LEMMA 1.

(15) $$\tilde{s}_{(j)}(\Lambda) = \begin{vmatrix} 1 & 0 & \cdots & 1 & 0 \\ \lambda_1 & 1 & \cdots & \lambda_N & 1 \\ \lambda_1^2 & 2\lambda_1 & \cdots & \lambda_N^2 & 2\lambda_N \\ \vdots & \vdots & \cdots & \vdots & \vdots \\ \lambda_1^{2N-2} & (2N-2)\lambda_1^{2N-3} & \cdots & \lambda_N^{2N-2} & (2N-2)\lambda_N^{2N-3} \\ \lambda_1^{2N+j-1} & (2N+j-1)\lambda_1^{2N+j-2} & \cdots & \lambda_N^{2N+j-1} & (2N+j-1)\lambda_N^{2N+j-2} \end{vmatrix}$$
$$\times V(\lambda)^{-4},$$

with the $(k, 2j-1)$ entry of the matrix being a power of $\lambda_j$ with the exponent $k-1$ if $k \neq 2N$ and $2N+j-1$ if $k = 2N$, and the $(k, 2j)$ entry being the derivative of the $(k, 2j-1)$ entry with respect to $\lambda_j$.

To prove this lemma, we need the well-known fact (see [23]), proven by L'Hôpital's rule

(16) $$V(\lambda)^4 = \begin{vmatrix} 1 & 0 & \cdots & 1 & 0 \\ \lambda_1 & 1 & \cdots & \lambda_N & 1 \\ \vdots & \vdots & \cdots & \vdots & \vdots \\ \lambda_1^{2N-1} & (2N-1)\lambda_1^{2N-2} & \cdots & \lambda_N^{2N-1} & (2N-1)\lambda_N^{2N-2} \end{vmatrix},$$

with the $(k, 2j-1)$ entry being $\lambda_j^{k-1}$ and the $(k, 2j)$ entry $(k-1)\lambda_j^{k-2}$.

PROOF OF LEMMA 1. Applying the L'Hôpital's rule repeatedly with respect to $x_{2i}$, $i = 1, \ldots, N$, we get the identity

$$\begin{vmatrix} 1 & 0 & \cdots & 1 & 0 \\ \lambda_1 & 1 & \cdots & \lambda_N & 1 \\ \lambda_1^2 & 2\lambda_1 & \cdots & \lambda_N^2 & 2\lambda_N \\ \vdots & \vdots & \cdots & \vdots & \vdots \\ \lambda_1^{2N-2} & (2N-2)\lambda_1^{2N-3} & \cdots & \lambda_N^{2N-2} & (2N-2)\lambda_N^{2N-3} \\ \lambda_1^{2N+j-1} & (2N+j-1)\lambda_1^{2N+j-2} & \cdots & \lambda_N^{2N+j-1} & (2N+j-1)\lambda_N^{2N+j-2} \end{vmatrix}$$

$$\Big/ \begin{vmatrix} 1 & 0 & \cdots & 1 & 0 \\ \lambda_1 & 1 & \cdots & \lambda_N & 1 \\ \vdots & \vdots & \cdots & \vdots & \vdots \\ \lambda_1^{2N-1} & (2N-1)\lambda^{2N-2} & \cdots & \lambda_N^{2N-1} & (2N-1)\lambda_N^{2N-2} \end{vmatrix}$$

$$= \left( \frac{\partial^N}{\partial x_2 \partial x_4 \cdots \partial x_{2N}} \begin{vmatrix} 1 & 1 & \cdots & 1 & 1 \\ x_1 & x_2 & \cdots & x_{2N-1} & x_{2N} \\ \vdots & \vdots & \cdots & \vdots & \vdots \\ x_1^{2N-2} & x_2^{2N-2} & \cdots & x_{2N-1}^{2N-2} & x_{2N}^{2N-2} \\ x_1^{2N+j-1} & x_2^{2N+j-1} & \cdots & x_{2N-1}^{2N+j-1} & x_{2N}^{2N+j-1} \end{vmatrix} \right.$$

$$\left. \Big/ \frac{\partial^N}{\partial x_2 \partial x_4 \cdots \partial x_{2N}} \begin{vmatrix} 1 & 1 & \cdots & 1 & 1 \\ x_1 & x_2 & \cdots & x_{2N-1} & x_{2N} \\ \vdots & \vdots & \cdots & \vdots & \vdots \\ x_1^{2N-2} & x_2^{2N-2} & \cdots & x_{2N-1}^{2N-2} & x_{2N}^{2N-2} \\ x_1^{2N+j-1} & x_2^{2N+j-1} & \cdots & x_{2N-1}^{2N+j-1} & x_{2N}^{2N+j-1} \end{vmatrix} \Bigg|_{\substack{x_{2i-1}=x_{2i}=\lambda_i \\ i=1,\ldots,N}} \right)$$

$$= s_{(j)}(\lambda_1, \lambda_1, \lambda_2, \lambda_2, \ldots, \lambda_N, \lambda_N) = \tilde{s}_{(j)}(\Lambda),$$

from the matrix representation of Schur polynomials, and now use (16) to get the compact formula (15). □

Substituting (15) into (14), we get

$$V(\lambda)^4 \int_{Q \in Sp(N)} e^{2M\Re \operatorname{Tr}((I-\Sigma^{-1})Q\Lambda Q^{-1})} \, dQ$$

(17) $$= \begin{vmatrix} 1 & 0 & \cdots & 1 & 0 \\ \lambda_1 & 1 & \cdots & \lambda_N & 1 \\ \lambda_1^2 & 2\lambda_1 & \cdots & \lambda_N^2 & 2\lambda_N \\ \vdots & \vdots & \cdots & \vdots & \vdots \\ \lambda_1^{2N-2} & (2N-2)\lambda_1^{2N-3} & \cdots & \lambda_N^{2N-2} & (2N-2)\lambda_N^{2N-3} \\ p(\lambda_1) & p'(\lambda_1) & \cdots & p(\lambda_N) & p'(\lambda_N) \end{vmatrix}$$

$$= \frac{1}{C} \begin{vmatrix} 1 & 0 & \cdots & 1 & 0 \\ \lambda_1 & 1 & \cdots & \lambda_N & 1 \\ \lambda_1^2 & 2\lambda_1 & \cdots & \lambda_N^2 & 2\lambda_N \\ \vdots & \vdots & \cdots & \vdots & \vdots \\ \lambda_1^{2N-2} & (2N-2)\lambda_1^{2N-3} & \cdots & \lambda_N^{2N-2} & (2N-2)\lambda_N^{2N-3} \\ e^{a/(1+a)2M\lambda_1} & \frac{a}{1+a}2Me^{a/(1+a)2M\lambda_1} & \cdots & e^{a/(1+a)2M\lambda_N} & \frac{a}{1+a}2Me^{a/(1+a)2M\lambda_N} \end{vmatrix}$$

$$= \frac{1}{C} \tilde{V}^4(\lambda),$$

where

$$p(x) = \sum_{j=0}^{\infty} \frac{1}{\prod_{i=0}^{j-1}(2N+i)} \left(\frac{a}{1+a}2M\right)^j x^{2N+j-1}$$

$$= \frac{(2N-1)!}{(a/(1+a)2M)^{2N-1}} \left( e^{a/(1+a)2Mx} - \sum_{j=0}^{2N-2} \frac{1}{j!} \left(\frac{a}{1+a}2Mx\right)^j \right),$$

and if $k \neq 2N$, the $(k, 2j-1)$ entries in both matrices are $\lambda_j^{k-1}$, and the $(k, 2j)$ entries are $(k-1)\lambda_j^{k-2}$, and the $2N, 2i-1$ entry in the former (latter) matrix is $p(\lambda_i)$ (resp. $e^{a/(1+a)2M\lambda_i}$) and the $2N, 2i$ entry $p'(\lambda_i)$ (resp. $\frac{a}{1+a}2Me^{a/(1+a)2M\lambda_i}$).

PROOF OF THEOREM 3. Formulas (8), (10) and (17) together give the result (7). □

2.2. *The Pfaffian and determinantal formulas.* With the formula (7) ready to use, we apply the standard RMT technique to get the distribution formula for the largest sample eigenvalue, in the same spirit as the solution of the LSE. Our process below is closely parallel to that in [31] to the LSE.

First, we find a skew orthogonal basis $\{\varphi_0(x), \varphi_1(x), \ldots, \varphi_{2N-1}(x)\}$ of the linear space spanned by $\{1, x, x^2, \ldots, x^{2N-2}, e^{a/(1+a)2Mx}\}$. We require that the $\varphi_j(x)$ is a linear combination of $\{1, x, x^2, \ldots, x^j\}$ if $j < 2N-1$, while $\varphi_{2N-1}(x)$ can be arbitrary, with the skew inner products among them

$$\langle \varphi_j(x), \varphi_k(x) \rangle_4 = \int_0^{\infty} (\varphi_j(x)\varphi_k'(x) - \varphi_j'(x)\varphi_k(x)) x^{2(M-N)+1} e^{-2Mx} dx$$

$$= \begin{cases} r_{j/2}, & \text{if } j \text{ is even and } k = j+1, \\ -r_{k/2}, & \text{if } k \text{ is even and } j = k+1, \\ 0, & \text{otherwise.} \end{cases}$$

Then we can reformulate the distribution function of $\lambda$ as

$$P(\lambda) = \frac{1}{C} \begin{vmatrix} \varphi_0(\lambda_1) & \varphi_0'(\lambda_1) & \cdots & \varphi_0(\lambda_N) & \varphi_0'(\lambda_N) \\ \varphi_1(\lambda_1) & \varphi_1'(\lambda_1) & \cdots & \varphi_1(\lambda_N) & \varphi_1'(\lambda_N) \\ \vdots & \vdots & \cdots & \vdots & \vdots \\ \varphi_{2N-1}(\lambda_1) & \varphi_{2N-1}'(\lambda_1) & \cdots & \varphi_{2N-1}(\lambda_N) & \varphi_{2N-1}'(\lambda_N) \end{vmatrix}$$

(18)
$$\times \prod_{j=1}^{N} (\lambda_j^{2(M-N)+1} e^{-2M\lambda_j})$$

$$= \frac{1}{C} \begin{vmatrix} \psi_0(\lambda_1) & \psi_0'(\lambda_1) & \cdots & \psi_0(\lambda_N) & \psi_0'(\lambda_N) \\ \psi_1(\lambda_1) & \psi_1'(\lambda_1) & \cdots & \psi_1(\lambda_N) & \psi_1'(\lambda_N) \\ \vdots & \vdots & \cdots & \vdots & \vdots \\ \psi_{2N-1}(\lambda_1) & \psi_{2N-1}'(\lambda_1) & \cdots & \psi_{2N-1}(\lambda_N) & \psi_{2N-1}'(\lambda_N) \end{vmatrix},$$

where

(19) $$\psi_i(x) = \varphi_i(x) x^{M-N+1/2} e^{-Mx}.$$

For an arbitrary function $f(x)$ on $[0, \infty)$, by the formula of de Bruijn [9],

$$\int_0^\infty \cdots \int_0^\infty \begin{vmatrix} \psi_0(\lambda_1) & \psi_0'(\lambda_1) & \cdots & \psi_0(\lambda_N) & \psi_0'(\lambda_N) \\ \psi_1(\lambda_1) & \psi_1'(\lambda_1) & \cdots & \psi_1(\lambda_N) & \psi_1'(\lambda_N) \\ \vdots & \vdots & \cdots & \vdots & \vdots \\ \psi_{2N-1}(\lambda_1) & \psi_{2N-1}'(\lambda_1) & \cdots & \psi_{2N-1}(\lambda_N) & \psi_{2N-1}'(\lambda_N) \end{vmatrix}$$

(20)
$$\times \prod_{i=1}^{N} (1 + f(\lambda_i)) \, d\lambda_i = C \operatorname{Pf}(P(1+f)),$$

where $P(1+f)$ is a $2N \times 2N$ matrix, whose entries depend on $1+f$ in the following way:

$$(P(1+f))_{j,k} = \int_0^\infty (\psi_{j-1}(x)\psi_{k-1}'(x) - \psi_{j-1}'(x)\psi_{k-1}(x))(1+f(x)) \, dx.$$

Now we define a matrix $Z$ as

$$Z = \begin{pmatrix} 0 & r_0 & & & & \\ -r_0 & 0 & & & & \\ & & 0 & r_1 & & \\ & & -r_1 & 0 & & \\ & & & & \ddots & \\ & & & & & 0 & r_{N-1} \\ & & & & & -r_{N-1} & 0 \end{pmatrix},$$

with

$$Z_{j,k} = \begin{cases} r_{k/2-1}, & \text{if } k \text{ is even and } j = k-1, \\ -r_{j/2-1}, & \text{if } j \text{ is even and } k = j-1, \\ 0, & \text{otherwise}, \end{cases}$$

and define for $j = 0, \ldots, N-1$, $\eta = Z^{-1}\psi$, that is,
$$\eta_{2j}(x) = -\frac{\psi_{2j+1}(x)}{r_j} \quad \text{and} \quad \eta_{2j+1}(x) = \frac{\psi_{2j}(x)}{r_j}.$$

So we have
$$\begin{aligned}(P(1+f))_{j,k} &= \int_0^\infty (\psi_{j-1}(x)\psi'_{k-1}(x) - \psi'_{j-1}(x)\psi_{k-1}(x))\,dx \\ &\quad + \int_0^\infty (\psi_{j-1}(x)\psi'_{k-1}(x) - \psi'_{j-1}(x)\psi_{k-1}(x))f(x)\,dx \\ &= Z_{j,k} + \int_0^\infty (\psi_{j-1}(x)\psi'_{k-1}(x) - \psi'_{j-1}(x)\psi_{k-1}(x))f(x)\,dx.\end{aligned}$$

And if we denote $Q(1+f) = Z^{-1}P(1+f)$, then
$$Q(1+f)_{j,k} = \delta_{j,k} + \int_0^\infty (\eta_{j-1}(x)\psi'_{k-1}(x) - \eta'_{j-1}(x)\psi_{k-1}(x))f(x)\,dx.$$

If we choose $f$ to be $-\chi_{(T,\infty)}$, then the integral on the left-hand side of (20), after multiplying a constant, is the probability of all $\lambda_i$'s smaller than $T$. In latter part of the paper, we abbreviate $\chi_{(T,\infty)}$ to $\chi$. So we get for a $T$-independent constant
$$\mathbb{P}(\max(\lambda_i) \leq T) = C\operatorname{Pf}(P(1-\chi)),$$
and
$$(\mathbb{P}(\max(\lambda_i) \leq T))^2 = C^2 \det(P(1-\chi)) = C^2 \det(Q(1-\chi)).$$

In linear algebra, we have the determinant identity
$$\det(I - AB) = \det(I - BA), \tag{21}$$

for $A$ an $m \times n$ matrix and $B$ an $n \times m$ matrix, but the identity still holds in infinite dimensional settings [16]. Letting det mean a Fredholm determinant for matrix integral operators, we describe a setting due to Tracy–Widom [31].

If $A$ is an operator from $L^2([0,\infty)) \times L^2([0,\infty))$ to the vector space $\mathbb{R}^{2N}$ with
$$A\begin{pmatrix} g(x) \\ h(x) \end{pmatrix}_j = \int_0^\infty \chi(x)\eta_{j-1}(x)g(x)\,dx - \int_0^\infty \chi(x)\eta'_{j-1}(x)h(x)\,dx,$$

and $B$ is an operator from $\mathbb{R}^{2N}$ to $L^2([0,\infty)) \times L^2([0,\infty))$ with
$$B\begin{pmatrix} c_1 \\ \vdots \\ c_{2N} \end{pmatrix} = \begin{pmatrix} \sum_{k=1}^{2N} c_k \psi'_{k-1}(x)\chi(x) \\ \sum_{k=1}^{2N} c_k \psi_{k-1}(x)\chi(x) \end{pmatrix},$$

then
$$I - AB = Q(1 - \chi)$$

and
$$I - BA = I - \chi(x) \begin{pmatrix} S_4(x,y) & SD_4(x,y) \\ IS_4(x,y) & S_4(y,x) \end{pmatrix} \chi(y),$$

where $S_4(x,y)$, $IS_4(x,y)$ and $SD_4(x,y)$ are integral operators whose kernels are

$$S_4(x,y) = \sum_{j=0}^{2N-1} \psi'_j(x) \eta_j(y)$$
$$= \sum_{j=0}^{N-1} \frac{1}{r_j}(-\psi'_{2j}(x)\psi_{2j+1}(y) + \psi'_{2j+1}(x)\psi_{2j}(y)),$$

$$SD_4(x,y) = \sum_{j=0}^{2N-1} -\psi'_j(x)\eta'_j(y)$$

(22)
$$= \sum_{j=0}^{N-1} \frac{1}{r_j}(\psi'_{2j}(x)\psi'_{2j+1}(y) - \psi'_{2j+1}(x)\psi'_{2j}(y)),$$

$$IS_4(x,y) = \sum_{j=0}^{2N-1} \psi_j(x)\eta_j(y)$$
$$= \sum_{j=0}^{N-1} \frac{1}{r_j}(-\psi_{2j}(x)\psi_{2j+1}(y) + \psi_{2j+1}(x)\psi_{2j}(y)),$$

$$S_4(y,x) = \sum_{j=0}^{2N-1} -\psi_j(x)\eta'_j(y)$$

(23)
$$= \sum_{j=0}^{N-1} \frac{1}{r_j}(\psi_{2j}(x)\psi'_{2j+1}(y) - \psi_{2j+1}(x)\psi'_{2j}(y)).$$

REMARK 5. It is clear that the nomenclature of $SD_4(x,y)$ is due to the fact that $SD_4(x,y)$ is the negative of the derivative of $S_4(x,y)$. But $IS_4(x,y)$, which gets its name in the same way in earlier literature in GSE (e.g., [30]), in our problem may not satisfy the equation

$$IS_4(x,y) = -\int_x^\infty S_4(t,y)\,dt,$$

since the integral on the right-hand side may diverge.

In conclusion,

$$(\mathbb{P}(\max(\lambda_i) \leq T))^2 = C^2 \det\left(I - \chi(x)\begin{pmatrix} S_4(x,y) & SD_4(x,y) \\ IS_4(x,y) & S_4(y,x) \end{pmatrix}\chi(y)\right),$$

and we can find that $C^2 = 1$ by taking the limit $T \to \infty$. We define a $2 \times 2$ matrix kernel as

$$\begin{aligned} P_T(x,y) &= \chi(x)\begin{pmatrix} S_4(x,y) & SD_4(x,y) \\ IS_4(x,y) & S_4(y,x) \end{pmatrix}\chi(y) \\ &= \begin{pmatrix} \chi(x)S_4(x,y)\chi(y) & \chi(x)DS_4(x,y)\chi(y) \\ \chi(x)IS_4(x,y)\chi(y) & \chi(x)S_4(y,x)\chi(y) \end{pmatrix}, \end{aligned}$$

then we have

$$(\mathbb{P}(\max(\lambda_i) \leq T))^2 = \det(I - P_T(x,y)).$$

2.3. $S_4(x,y)$ *in terms of Laguerre polynomials.* In manipulation of skew orthogonal polynomials, we take the approach of [1], and all classical orthogonal polynomial properties are from [28].

Since Laguerre polynomials by definition satisfy the orthogonal property

$$\int_0^\infty L_j^{(\alpha)} L_k^{(\alpha)} x^\alpha e^{-x}\, dx = \frac{(j+\alpha)!}{j!}\delta_{j,k},$$

and they have the differential identity [we assume $L_n^{(\alpha)}(x) = 0$ if $n < 0$]

(24) $$x\frac{d}{dx}L_n^{(\alpha)}(x) = nL_n^{(\alpha)}(x) - (n+\alpha)L_{n-1}^{(\alpha)}(x),$$

it is easy to get that

$$\begin{aligned} &\langle L_j^{(2(M-N))}(2Mx), L_k^{(2(M-N))}(2Mx)\rangle_4 \\ &= \int_0^\infty \bigg( L_j^{(2(M-N))}(2Mx)\frac{d}{dx}L_k^{(2(M-N))}(2Mx) \\ &\qquad - L_k^{(2(M-N))}(2Mx)\frac{d}{dx}L_j^{(2(M-N))}(2Mx)\bigg) \\ &\qquad \times x^{2(M-N)+1}e^{-2Mx}\, dx \\ &= \begin{cases} \left(\dfrac{1}{2M}\right)^{2(M-N)+1}\dfrac{(j+2(M-N))!}{(j-1)!}, & \text{if } j = k+1, \\ -\left(\dfrac{1}{2M}\right)^{2(M-N)+1}\dfrac{(k+2(M-N))!}{(k-1)!}, & \text{if } k = j+1, \\ 0, & \text{otherwise.} \end{cases} \end{aligned}$$

So we can choose for $j = 0, \ldots, N-2$,

$$\varphi_{2j}(x) = \sum_{k=0}^{j} \left( \prod_{i=1}^{k} \frac{2i-1}{2i+2(M-N)} \right) L_{2k}^{(2(M-N))}(2Mx), \tag{25}$$

$$\varphi_{2j+1}(x) = -L_{2j+1}^{(2(M-N))}(2Mx) \tag{26}$$

and

$$r_j = \left( \frac{1}{2M} \right)^{2(M-N)+1} \frac{(2j+2(M-N)+1)!}{(2j)!} \prod_{k=1}^{j} \frac{2k-1}{2k+2(M-N)}. \tag{27}$$

We can also choose

$$\varphi_{2N-2}(x) = \sum_{k=0}^{N-1} \left( \prod_{i=1}^{k} \frac{2i-1}{2i+2(M-N)} \right) L_{2k}^{(2(M-N))}(2Mx),$$

but $\varphi_{2N-1}(x)$ is not a polynomial and needs to be treated separately.

By the Rodrigues' representation

$$x^\alpha e^{-x} L_n^{(\alpha)}(x) = \frac{1}{n!} \frac{d^n}{dx^n}(e^{-x} x^{n+\alpha}),$$

and repeated integration by parts, we get for $n > 0$

$$\langle e^{a/(1+a)2Mx}, L_n^{(2(M-N))}(2Mx) \rangle_4$$
$$= \left( \frac{1+a}{2M} \right)^{2(M-N)+1} \left( (-a)^{n+1} \frac{(n+2(M-N)+1)!}{n!} \right.$$
$$\left. - (-a)^{n-1} \frac{(n+2(M-N))!}{(n-1)!} \right)$$

and

$$\langle e^{a/(1+a)2Mx}, L_0^{(2(M-N))}(2Mx) \rangle_4 = -\left( \frac{1+a}{2M} \right)^{2(M-N)+1} a(2(M-N)+1)!,$$

so that

$$\langle e^{a/(1+a)2Mx}, \varphi_{2j}(x) \rangle_4$$
$$= -\left( \frac{1+a}{2M} \right)^{2(M-N)+1} a^{2j+1} \frac{(2j+2(M-N)+1)!}{(2j)!} \prod_{k=1}^{j} \frac{2k-1}{2k+2(M-N)}$$

and

$$\langle e^{a/(1+a)2Mx}, \varphi_{2j+1}(x) \rangle_4$$
$$= -\left( \frac{1+a}{2M} \right)^{2(M-N)+1} \left( a^{2j+2} \frac{(2j+2(M-N)+2)!}{(2j+1)!} \right.$$
$$\left. - a^{2j} \frac{(2j+2(M-N)+1)!}{(2j)!} \right).$$

Now by the skew orthogonality, we can choose

$$\varphi_{2N-1}(x) = e^{a/(1+a)2Mx} - \sum_{j=0}^{N-2} \frac{1}{r_j}(\langle e^{a/(1+a)2Mx}, \varphi_{2j+1}(x)\rangle_4 \varphi_{2j}(x)$$
$$- \langle e^{a/(1+a)2Mx}, \varphi_{2j}(x)\rangle_4 \varphi_{2j+1}(x))$$
$$- (1+a)^{2(M-N)+1} a^{2N-2} \prod_{j=1}^{N-1} \frac{2j+2(M-N)}{2j-1} \varphi_{2N-2}(x)$$
$$= e^{a/(1+a)2Mx} - (1+a)^{2(M-N)+1} \sum_{j=0}^{2N-2} (-a)^j L_j^{(2(M-N))}(2Mx)$$

and

$$r_{N-1} = \left(\frac{1+a}{2M}\right)^{2(M-N)+1} a^{2N-1} \frac{(2M-1)!}{(2N-2)!} \prod_{k=1}^{N-1} \frac{2k-1}{2k+2(M-N)}.$$

Now, we write $S_4(x,y)$ as $S_{4a}(x,y) + S_{4b}(x,y)$, where

(28) $$S_{4a}(x,y) = \sum_{j=0}^{N-2} \frac{1}{r_j}(-\psi'_{2j}(x)\psi_{2j+1}(y) + \psi'_{2j+1}(x)\psi_{2j}(y))$$

and

(29) $$S_{4b}(x,y) = \frac{1}{r_{N-1}}(-\psi'_{2N-2}(x)\psi_{2N-1}(y) + \psi'_{2N-1}(x)\psi_{2N-2}(y)),$$

and simplify them separately.

The formula (28) of our $S_{4a}(x,y)$ is also the formula for $S_4(x,y)$ in the LSE problem, with parameters $M$ and $N-2$, and has been well studied. For completeness we derive its Laguerre polynomial expression here, following [1].

By the differential identity (24) and the identity

$$nL_n^{(\alpha)}(x) = (-x + 2n + \alpha - 1)L_{n-1}^{(\alpha)}(x) - (n + \alpha - 1)L_{n-2}^{(\alpha)}(x),$$

we get, remembering the definition (19), the telescoping sequence

$$\psi'_{2j}(x) = \sum_{k=0}^{j} \left(\prod_{i=1}^{k} \frac{2i-1}{2i+2(M-N)}\right.$$
$$\left. \times \left(M - N + 1/2 - Mx + x\frac{d}{dx}\right) L_{2k}^{(2(M-N))}(2Mx)\right)$$
$$\times x^{M-N-1/2} e^{-Mx}$$

$$= \frac{1}{2} \sum_{k=0}^{j} \prod_{i=1}^{k} \frac{2i-1}{2i+2(M-N)} ((2k+1)L_{2k+1}^{(2(M-N))}(2Mx)$$

(30)
$$- (2k+2(M-N))L_{2k-1}^{(2(M-N))}(2Mx))$$

$$\times x^{M-N-1/2}e^{-Mx}$$

$$= \frac{1}{2} \left( \prod_{k=1}^{j} \frac{2k-1}{2k+2(M-N)} \right)$$

$$\times (2j+1)L_{2j+1}^{(2(M-N))}(2Mx)x^{M-N-1/2}e^{-Mx}$$

and

$$\psi'_{2j+1}(x) = -\left( M - N + 1/2 - Mx + x\frac{d}{dx} \right) L_{2j+1}^{(2(M-N))}$$

$$\times (2Mx)x^{M-N-1/2}e^{-Mx}$$

(31)
$$= -\frac{1}{2}((2j+2)L_{2j+2}^{(2(M-N))}(2Mx)$$

$$- (2j+2(M-N)+1)L_{2j}^{(2(M-N))}(2Mx))$$

$$\times x^{M-N-1/2}e^{-Mx}.$$

Therefore, if we substitute (27), (30) and (31) into (28), we get after some trick,

$$S_{4a}(x,y) = \frac{1}{2}(2M)^{2(M-N)+1}x^{M-N-1/2}e^{-Mx}y^{M-N+1/2}e^{-My}$$

$$\times \left\{ \sum_{j=0}^{2N-2} \frac{j!}{(j+2(M-N))!}L_j^{(2(M-N))}(2Mx)L_j^{(2(M-N))}(2My) \right.$$

$$- \frac{(2N-2)!}{(2M-2)!} \left( \prod_{j=1}^{N-1} \frac{2j+2(M-N)}{2j-1} \right)$$

$$\left. \times L_{2N-2}^{(2(M-N))}(2Mx)\varphi_{2N-2}(y) \right\}.$$

Furthermore, we can simplify $\psi_{2N-2}(x)$. Since for $j \neq 2N - 1$ [if we define $\varphi_j(x)$ and then $\psi_j(x)$ for $j > 2N - 1$ by the formula (25) and (26)]

$$\int_0^\infty (\psi_{2N-2}(x)\psi'_j(x) - \psi'_{2N-2}(x)\psi_j(x))\,dx = 0,$$

we get for $j \neq 2N - 1$, using integration by parts,

$$\int_0^\infty \psi'_{2N-2}(x)L_j^{(2(M-N))}(2Mx)x^{M-N+1/2}e^{-Mx}\,dx = 0.$$

So by the orthogonal property of Laguerre polynomials, we get
$$\psi'_{2N-2}(x) = C L_{2N-1}^{(2(M-N))}(2Mx) x^{M-N-1/2} e^{-Mx},$$
and we can determine that
$$C = \frac{2N-1}{2} \prod_{j=1}^{N-1} \frac{2j-1}{2j+2(M-N)}$$
without much difficulty. Together with the fact $\lim_{x\to\infty} \psi_{2N-2}(x) = 0$, we get
$$\psi_{2N-2}(x) = -\frac{2N-1}{2} \prod_{j=1}^{N-1} \frac{2j-1}{2j+2(M-N)}$$
$$\times \int_x^\infty t^{M-N-1/2} e^{-Mt} L_{2N-1}^{(2(M-N))}(2Mt)\, dt.$$

Now, we can write $S_{4a}(x,y)$ as $S_{4a1}(x,y) + S_{4a2}(x,y)$, where
$$S_{4a1}(x,y) = \frac{1}{2}(2M)^{2(M-N)+1}$$
(32)
$$\times \sum_{j=0}^{2N-2} \frac{j!}{(j+2(M-N))!} L_j^{(2(M-N))}(2Mx) x^{M-N-1/2} e^{-Mx}$$
$$\times L_j^{(2(M-N))}(2My) y^{M-N+1/2} e^{-My}$$

and
$$S_{4a2}(x,y) = \frac{1}{4}(2M)^{2(M-N)+1} \frac{(2N-1)!}{(2M-2)!}$$
(33)
$$\times L_{2N-2}^{(2(M-N))}(2Mx) x^{M-N-1/2} e^{-Mx}$$
$$\times \int_y^\infty t^{M-N-1/2} e^{-Mt} L_{2N-1}^{(2(M-N))}(2Mt)\, dt.$$

Finally,
$$S_{4b}(x,y) = -\frac{1}{2}\left(\frac{2M}{1+a}\right)^{2(M-N)+1} a^{-(2N-1)} \frac{(2N-1)!}{(2M-1)!}$$
(34)
$$\times \left\{ L_{2N-1}^{(2(M-N))}(2Mx) x^{M-N-1/2} e^{-Mx} \psi_{2N-1}(y) \right.$$
$$\left. + \psi'_{2N-1}(x) \int_y^\infty L_{2N-1}^{(2(M-N))}(2Mt) t^{M-N-1/2} e^{-Mt}\, dt \right\},$$

and we can take the asymptotic analyses of $S_{4a1}(x,y)$, $S_{4a2}(x,y)$ and $S_{4b}(x,y)$ separately.

**3. Asymptotic analysis.** In order to consider the rescaled distribution problem, we wish to find the probability of the largest sample eigenvalue being in the domain $(0, p+qT]$. We can put the kernel in the new coordinate system [after a conjugation by $\begin{pmatrix} q^{1/2} & 0 \\ 0 & q^{-1/2} \end{pmatrix}$], and get

$$(\mathbb{P}(\max(\lambda_i) \le p + qT))^2 = \det\left(I - \chi(x)\begin{pmatrix} \widetilde{S}_4(\xi,\eta) & \widetilde{SD}_4(\xi,\eta) \\ \widetilde{IS}_4(\xi,\eta) & \widetilde{S}_4(\eta,\xi) \end{pmatrix}\chi(\eta)\right)$$
$$= \det(I - \widetilde{P}_T(\xi,\eta)),$$

where as $L^2$ functions,

(35) $$\widetilde{SD}_4(\xi,\eta) = q^2 SD_4(x,y)|_{\substack{x=p+q\xi \\ y=p+q\eta}},$$

(36) $$\widetilde{S}_4(\xi,\eta) = qS_4(x,y)|_{\substack{x=p+q\xi \\ y=p+q\eta}},$$

(37) $$\widetilde{IS}_4(\xi\eta) = IS_4(x,y)|_{\substack{x=p+q\xi \\ y=p+q\eta}}$$

and

$$\widetilde{P}_T(\xi,\eta) = \chi(\xi)\begin{pmatrix} \widetilde{S}_4(\xi,\eta) & \widetilde{SD}_4(\xi,\eta) \\ \widetilde{IS}_4(\xi,\eta) & \widetilde{S}_4(\eta,\xi) \end{pmatrix}\chi(\eta).$$

In this section, we want to prove that for fixed $\gamma \ge 1$ and $a > -1$, we can choose suitable $p_M$ and $q_M$ depending on $M$, so that for any $T$,

$$\lim_{M\to\infty}(\mathbb{P}(\max(\lambda_i) \le p_M + q_M T))^2 = \lim_{M\to\infty}\det(I - \widetilde{P}_T(\xi,\eta)) = f_a(T),$$

where $f_a$ is a function to be determined.

To prove the convergence of Fredholm determinants, we may use that $\widetilde{P}_T(\xi,\eta)$ is in trace class for any $M$ and converges to a certain $2 \times 2$ matrix kernel in trace norm. Equivalently, we may use that each entry of $\widetilde{P}_T(\xi,\eta)$ is in trace class and converges to a scalar kernel in trace norm. It turns out later that the $\widetilde{P}_T(\xi,\eta)$'s may not satisfy these requirements, but certain conjugates do.

Since the $IS_4(x,y)$ and $DS_4(x,y)$ are of the same form as $S_4(x,y)$, we only show the asymptotic analysis of $S_4(x,y)$, and state the result for the other two, for which the arguments are the same.

3.1. *Proof of the $-1 < a < \gamma^{-1}$ part of Theorem 1.* In case $-1 < a \le \gamma^{-1}$, we choose $p_M = (1+\gamma^{-1})^2$ and $q_M = \frac{(1+\gamma)^{4/3}}{\gamma(2M)^{2/3}}$, and denote [here $*$ stands for 4, 4a, 4a1, 4a2 and 4b; the definition of $\widetilde{S}_*(\xi,\eta)$ in (38) is only used in

Sections 3.1 and 3.2]

$$\widetilde{S}_*(\xi,\eta) = \frac{(1+\gamma)^{4/3}}{\gamma(2M)^{2/3}} S_*(x,y)|_{\substack{x=(1+\gamma^{-1})^2+(1+\gamma)^{4/3}/(\gamma(2M)^{2/3})\xi \\ y=(1+\gamma^{-1})^2+(1+\gamma)^{4/3}/(\gamma(2M)^{2/3})\eta}} \tag{38}$$

$S_{4a}(x,y)$ is the formula for the upper-left entry of the $2 \times 2$ matrix kernel of the LSE problem with parameters $M$ and $N-1$, and its asymptotic behavior is well studied [14]. We want to prove that as $M \to \infty$, $S_{4a}(x,y)$ dominates $S_4(x,y)$ in the domain that we are interested in, and so naturally the distribution of the largest sample eigenvalue in the perturbed problem is the same as that in the LSE problem. (The difference between $N$ and $N-1$ is negligible.)

$S_{4a1}(x,y)$ is almost the kernel for the LUE problem with parameters $2M-2$ and $2N-2$, besides a factor $\sqrt{y/x}/2$. From a standard result for LUE [11], $\chi_T(\xi)\widetilde{S}_{4a1}(\xi,\eta)\chi_T(\eta)$ is in trace class and converges in trace norm to half of the Airy kernel

$$\lim_{M\to\infty} \chi(\xi)\widetilde{S}_{4a1}(\xi,\eta)\chi(\eta) = \frac{1}{2}\chi(\xi)K_{\mathrm{Airy}}(\xi,\eta)\chi(\eta). \tag{39}$$

More discussion see the Appendix.

For the $S_{4a2}(x,y)$ part, we also have in trace norm [14],

$$\lim_{M\to\infty} \chi(\xi)\widetilde{S}_{4a2}(\xi,\eta)\chi(\eta) = -\frac{1}{4}\chi(\xi)\operatorname{Ai}(\xi)\int_\eta^\infty \operatorname{Ai}(t)\,dt\,\chi(\eta). \tag{40}$$

We just sketch the proof. Since $\widetilde{S}_{4a2}(\xi,\eta)$ is a rank 1 operator, for the trace norm convergence, we only need to prove that in $L^2$ norm as functions in $\xi$ and respectively $\eta$,

$$\lim_{M\to\infty} \gamma^{-2N}(1+\gamma)^{4/3}(2M)^{1/3}e^{M-N} \\ \times L_{2N-2}^{(2(M-N))}(2Mx)x^{M-N-1/2}e^{-Mx}\chi(\xi) \\ = \operatorname{Ai}(\xi)\chi(\xi), \tag{41}$$

$$\lim_{M\to\infty} \gamma^{-2N}2Me^{M-N}\int_y^\infty L_{2N-1}^{(2(M-N))}(2Mt)t^{M-N-1/2}e^{-Mt}\,dt\,\chi(\eta) \\ = -\int_\eta^\infty \operatorname{Ai}(t)\,dt\,\chi(\eta) \tag{42}$$

and by the Stirling's formula,

$$\lim_{M\to\infty} (2M)^{2(M-N)-1}\frac{(2N-1)!}{(2M-2)!}e^{2(N-M)}\gamma^{4N-1} = 1.$$

By (33), (38), (41) and (42), we get

$$\chi(\xi)\widetilde{S}_{4a2}(\xi,\eta)\chi(\eta) = \frac{1}{4}(2M)^{2(M-N)-1}\frac{(2N-1)!}{(2M-2)!}e^{2(N-M)}\gamma^{4N-1}\gamma^{-2N}$$
$$\times (1+\gamma)^{4/3}(2M)^{1/3}e^{M-N}L_{2N-2}^{(2(M-N))}(2Mx)$$
$$\times x^{M-N-1/2}e^{-Mx}\chi(\xi)\gamma^{-2N}2Me^{M-N}$$
$$\times \int_y^\infty L_{2N-1}^{(2(M-N))}(2Mt)t^{M-N-1/2}e^{-Mt}\,dt\,\chi(\eta).$$

Therefore we get the trace norm convergence from the $L^2$ convergence by the fact that if $f_n(x) \to f(x)$ and $g_n(y) \to g(y)$ in $L^2$ norm, then we have the convergence of integral operators in trace norm:

$$f_n(x)g_n(y) \to f(x)g(y).$$

Finally, we need to analyze the term $S_{4b}(\xi,\eta)$, new to the perturbed problem. We need the following results:

PROPOSITION 5. *For fixed $\gamma \geq 1$ and $-1 < a < \gamma^{-1}$ and any $T$, we have the convergences in $L^2$ norm with respect to $\xi$ or $\eta$:*

$$\lim_{M\to\infty} \gamma^{-2N-1}(1+\gamma)^{4/3}(2M)^{1/3}e^{M-N}$$
$$\times L_{2N-1}^{(2(M-N))}(2Mx)x^{M-N-1/2}e^{-Mx}\chi(\xi)$$
$$= -\operatorname{Ai}(\xi)\chi(\xi),$$

$$\lim_{M\to\infty} \gamma^{-2N}2Me^{M-N}\int_y^\infty L_{2N-1}^{(2(M-N))}(2Mt)t^{M-N-1/2}e^{-Mt}\,dt\,\chi(\eta)$$
$$= -\int_\eta^\infty \operatorname{Ai}(t)\,dt\,\chi(\eta),$$

(43)
$$\lim_{M\to\infty} (1+a)^{2(N-M)-1}a^{-2N+1}\frac{(1-a\gamma)(2M)^{1/3}}{(\gamma+1)^{2/3}\gamma^{2N-1}}e^{M-N}\psi_{2N-1}(y)\chi(\eta)$$
$$= \operatorname{Ai}(\eta)\chi(\eta),$$

$$\lim_{M\to\infty} (1+a)^{2(N-M)-1}a^{-2N+1}\frac{(1-a\gamma)(\gamma+1)^{2/3}}{\gamma^{2N}(2M)^{1/3}}e^{M-N}\psi'_{2N-1}(x)\chi(\xi)$$
$$= \operatorname{Ai}'(\xi)\chi(\xi).$$

PROOF. We just prove the identity (43), and others can be done in the same way.

By the integral representation of Laguerre polynomials,

$$L_n^{(2(M-N))}(2Mx) = \frac{1}{2\pi i} \oint_C \frac{e^{-2Mxz}(z+1)^{n+2(M-N)}}{z^{n+1}} dz, \tag{44}$$

where $C$ is a contour around the pole 0, therefore we get

$$\begin{aligned}\varphi_{2N-1}(y) &= e^{a/(1+a)2My} \\ &\quad - \frac{(1+a)^{2(M-N)+1}}{2\pi i} \\ &\quad \times \oint_C e^{-2Myz} \frac{1+(a(z+1)/z)^{2N-1}}{1+a(z+1)/z} \frac{(z+1)^{2(M-N)}}{z} dz \\ &= e^{a/(1+a)2My} - \frac{(1+a)^{2(M-N)}}{2\pi i} \oint_C e^{-2Myz} \frac{(z+1)^{2(M-N)}}{z+a/(a+1)} dz \\ &\quad - \frac{(1+a)^{2(M-N)+1} a^{2N-1}}{2\pi i} \\ &\quad \times \oint_C e^{-2Myz} \frac{(z+1)^{2M}}{z^{2N}} \frac{z}{((a+1)z+a)(z+1)} dz.\end{aligned} \tag{45}$$

If the pole $z = -\frac{a}{a+1}$ is inside of $C$, then

$$\frac{(1+a)^{2(M-N)}}{2\pi i} \oint_C e^{-2Myz} \frac{(z+1)^{2(M-N)}}{z+a/(a+1)} dz = e^{a/(a+1)2My}$$

and

$$\begin{aligned}\varphi_{2N-1}(y) &= -\frac{(1+a)^{2(M-N)+1} a^{2N-1}}{2\pi i} \\ &\quad \times \oint_C e^{-2Myz} \frac{(z+1)^{2M}}{z^{2N}} \frac{z}{((a+1)z+a)(z+1)} dz.\end{aligned} \tag{46}$$

In later part of the proof, we make this condition hold, and will not mention the canceled terms, and we are then free to deform $C$ in (46) as we wish, provided it includes 0. We then proceed to a stationary phase analysis.

Since

$$e^{-2Myz} \frac{(z+1)^{2M}}{z^{2N}} = e^{2M(-(1+\gamma^{-1})^2 z + \log(z+1) - \gamma^{-2} \log z) - \frac{(1+\gamma)^{4/3}}{\gamma}(2M)^{1/3}\eta z}$$

(we do not need to concern ourselves about the ambiguity of the value of logarithmic functions), if we denote

$$f(z) = -(1+\gamma^{-1})^2 z + \log(z+1) - \gamma^{-2} \log z, \tag{47}$$

then we get:

- $f'(z) = -\frac{((1+\gamma^{-1})z+\gamma^{-1})}{z(z+1)}$, with the zero point $z = -\frac{1}{\gamma+1}$;
- $f''(-\frac{1}{\gamma+1}) = 0$;
- $f'''(-\frac{1}{\gamma+1}) = \frac{2(\gamma+1)^4}{\gamma^3} > 0$.

So locally around $z = -\frac{1}{\gamma+1}$,

$$f\left(-\frac{1}{\gamma+1} + w\right) = \frac{\gamma+1}{\gamma^2} + \log\gamma - (1-\gamma^{-2})\log(\gamma+1) + \gamma^{-2}\pi i$$
$$(48) \qquad\qquad + \frac{(\gamma+1)^4}{3\gamma^3} w^3 + R_1(w),$$

where

$$(49) \qquad\qquad R_1(w) = O(w^4), \qquad \text{as } w \to 0.$$

After the substitution $w = z + \frac{1}{\gamma+1}$, we get

$$\oint_C e^{2M(-yz + \log(z-1) - \gamma^{-2}\log z)} \frac{z}{((a+1)z + a)(z+1)} dz$$

$$= \oint_{\Gamma^M} \exp\Big\{ 2M\Big(\frac{\gamma+1}{\gamma^2} + \log\gamma - (1-\gamma^{-2})\log(\gamma+1) + \gamma^{-2}\pi i$$
$$\qquad\qquad + \frac{(\gamma+1)^4}{3\gamma^3} w^3 + R_1(w) - \frac{(1+\gamma)^{4/3}}{\gamma(2M)^{2/3}}\eta\Big(w - \frac{1}{\gamma+1}\Big)\Big)\Big\}$$
$$\qquad\qquad \times \frac{w - 1/(\gamma+1)}{((a+1)w + (a\gamma-1)/(\gamma+1))(w + \gamma/(\gamma+1))} dw$$

$$= -\frac{1}{a+1} \frac{\gamma^{2M-1}}{(\gamma+1)^{2(M-N)}} e^{2M/(1+\gamma)y}$$

$$\qquad \times \oint_{\Gamma^M} \exp\Big\{ \frac{-(1+\gamma)^{4/3}}{\gamma}(2M)^{1/3}\eta w + \frac{(1+\gamma)^4}{3\gamma^3} 2Mw^3 + 2MR_1(w)\Big\}$$

$$\qquad\qquad \times \frac{-(\gamma+1)w + 1}{(\gamma+1)/\gamma w + 1} \frac{1}{w + (a\gamma-1)/((\gamma+1)(a+1))} dw,$$

where $\Gamma^M$ is a contour around $\frac{1}{\gamma+1}$, composed of $\Gamma_1^M$, $\Gamma_2^M$, $\Gamma_3^M$ and $\Gamma_4^M$, which are defined as (see Figure 3)

$$\Gamma_1^M = \Big\{ (4-t)\frac{\gamma}{\gamma+1} e^{\pi i/3} \Big| 0 \le t \le 4 - \frac{1}{(1+\gamma)^{1/3}}(2M)^{-1/3} \Big\},$$

$$\Gamma_2^M = \Big\{ \frac{\gamma}{(1+\gamma)^{4/3}}(2M)^{-1/3} e^{-t\pi i} \Big| -\frac{1}{3} \le t \le \frac{1}{3} \Big\},$$

$$\Gamma_3^M = \Big\{ t\frac{\gamma}{\gamma+1} e^{5\pi i/3} \Big| \frac{1}{(1+\gamma)^{1/3}}(2M)^{-1/3} \le t \le 4 \Big\},$$

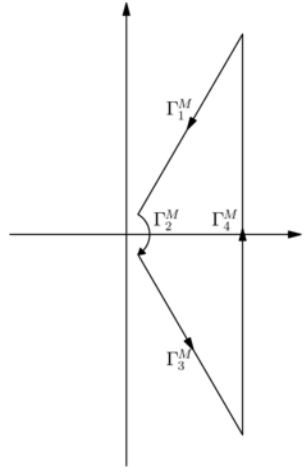

Fig. 3.  $\Gamma^M$.

$$\Gamma_4^M = \left\{2\frac{\gamma}{\gamma+1} + it \Big| -2\sqrt{3}\frac{\gamma}{\gamma+1} \leq t \leq 2\sqrt{3}\frac{\gamma}{\gamma+1}\right\}.$$

For the asymptotic analysis, we define

(50) $\quad \Gamma_{\text{local}}^M = \{z \in \Gamma^M | \Re(z) \leq (2M)^{-10/39}\}, \quad \Gamma_{\text{remote}}^M = (\Gamma_1^M \cup \Gamma_3^M) \setminus \Gamma_{\text{local}}^M,$

(51) $\quad\quad \Gamma_{<c}^\infty = \{w \in \Gamma^\infty | \Re(w) < c\}, \quad \Gamma_{\geq c}^\infty = \Gamma^\infty \setminus \Gamma_{<c}^\infty.$

Now, we denote

$$F_{aM}(\eta, w) = \frac{(1+\gamma)^{4/3}}{\gamma(2M)^{1/3}}$$

$$\times \exp\left\{-\frac{(1+\gamma)^{4/3}}{\gamma}(2M)^{1/3}\eta w + \frac{(1+\gamma)^4}{3\gamma^3}2Mw^3 + 2MR_1(w)\right\}$$

$$\times \frac{-(\gamma+1)w+1}{(\gamma+1)/\gamma w+1} \frac{1}{w + (a\gamma-1)/((\gamma+1)(a+1))},$$

and establish several lemmas for the proof. □

LEMMA 2. *If $T$ is fixed and $M$ is large enough, then for any $\eta > T$,*

$$\left|\frac{1}{2\pi i}\int_{\Gamma_4^M} F_{aM}(\eta, w)\, dw\right| < \frac{1}{3}\frac{e^{-\eta/2}}{M^{1/40}}.$$

PROOF. By (48) and (47),

$$\frac{(1+\gamma)^4}{3\gamma^3}2Mw^3 + 2MR_1(w)$$

$$= 2M\left(f\left(-\frac{1}{\gamma+1} + w\right) - \frac{\gamma+1}{\gamma^2}\right.$$

(52)
$$\left. - \log\gamma + (1 - \gamma^{-2})\log(\gamma+1) - \gamma^{-2}\pi i\right)$$

$$= 2M\left(-\frac{(\gamma+1)^2}{\gamma^2}w + \log\left(\frac{\gamma+1}{\gamma}w + 1\right)\right.$$

$$\left. - \gamma^{-2}\log((\gamma+1)w - 1) - \gamma^{-2}\pi i\right).$$

If $w \in \Gamma_4^M$, $\Re(w) = 2\frac{\gamma}{\gamma+1}$, and denote $\theta = \arg(w) \in [-\frac{\pi}{3}, \frac{\pi}{3}]$, we have

$$\Re\left(\frac{(1+\gamma)^4}{3\gamma^3} 2Mw^3 + 2MR_1(w)\right)$$

$$= 2M\left(-2\frac{\gamma+1}{\gamma} + \log(\sqrt{3^2 + (2\tan\theta)^2})\right.$$

$$\left. - \gamma^{-2}\log(\sqrt{(1+2\gamma)^2 + (2\gamma\tan\theta)^2})\right)$$

$$\leq 2M\left(-2\frac{\gamma+1}{\gamma} + \log\sqrt{21} - \gamma^{-2}\log(1+2\gamma)\right) < (\log\sqrt{21} - 2)2M < 0.$$

So on $\Gamma_4^M$, if $\eta \geq T$, $0 < \varepsilon' < 2 - \log\sqrt{21}$ and $M$ large enough,

$$|F_{aM}(\eta, w)| < \frac{(1+\gamma)^{4/3}}{\gamma}(2M)^{1/3}$$

$$\times \exp\{-2(\eta - T)(1+\gamma)^{1/3}(2M)^{1/3}$$

$$+ ((\log\sqrt{21} - 2) - 2T(1+\gamma)^{-2/3})2M\}$$

$$\times \left|\frac{-(\gamma+1)w + 1}{(\gamma+1)/\gamma w + 1} \frac{1}{w + (a\gamma - 1)/((\gamma+1)(a+1))}\right|$$

$$< e^{-2(\eta-T)(1+\gamma)^{1/3}(2M)^{1/3}} e^{(\log\sqrt{21}-2+\varepsilon')2M},$$

where $\varepsilon'$ is a positive number and $\varepsilon' < 2 - \log\sqrt{21}$. If $M$ is large enough,

(53)
$$e^{(\log\sqrt{21}-2+\varepsilon')2M} < \frac{2\pi}{2\sqrt{3}\gamma/(\gamma+1)} \frac{1}{3} \frac{e^{-T/2}}{M^{1/40}},$$

$$e^{-2(\eta-T)(1+\gamma)^{1/3}(2M)^{1/3}} < e^{T/2}e^{-\eta/2},$$

and we get the result, since

(54) $$\left|\frac{1}{2\pi i}\int_{\Gamma_4^M} F_{aM}(\eta, w)\, dw\right| \leq \frac{2\sqrt{3}\gamma/(\gamma+1)}{2\pi} \max_{w \in \Gamma_4^M}|F_{aM}(\eta, w)|.$$ □

LEMMA 3. *If $T$ is fixed and $M$ is large enough, then for any $\eta > T$,*

$$\left| \frac{1}{2\pi i} \int_{\Gamma^M_{\text{remote}}} F_{aM}(\eta, w)\, dw \right| < \frac{1}{3} \frac{e^{-\eta/2}}{M^{1/40}}.$$

PROOF. For $w \in \Gamma^M_{\text{remote}}$, we denote $l = \Re(w) = \frac{|w|}{2}$. Since $\arg(w) = \pm\frac{\pi}{3}$, we get by (52)

$$\Re\left( \frac{(1+\gamma)^4}{3\gamma^3} 2Mw^3 + 2MR_1(w) \right)$$
$$= 2M\left( -\frac{(\gamma+1)^2}{\gamma^2} l + \frac{1}{2}\log\left(1 + 2\frac{\gamma+1}{\gamma}l + 4\left(\frac{\gamma+1}{\gamma}l\right)^2\right) \right.$$
$$\left. - \frac{\gamma^{-2}}{2}\log(1 - 2(\gamma+1)l + 4(\gamma+1)^2 l^2) \right).$$

Then we take derivative

$$\frac{d}{dl}\left( -\frac{(\gamma+1)^2}{\gamma^2} l + \frac{1}{2}\log\left(1 + 2\frac{\gamma+1}{\gamma}l + 4\left(\frac{\gamma+1}{\gamma}l\right)^2\right) \right.$$
$$\left. - \frac{\gamma^{-2}}{2}\log(1 - 2(\gamma+1)l + 4(\gamma+1)^2 l^2) \right)$$
$$= -8\frac{(\gamma+1)^4}{\gamma^3} l^2 \left( 1 - (\gamma-1)\frac{\gamma+1}{\gamma}l + 4\left(\frac{\gamma+1}{\gamma}l\right)^2 \right)$$
$$\times \left\{ \left(1 + 2\frac{\gamma+1}{\gamma}l + 4\left(\frac{\gamma+1}{\gamma}l\right)^2\right) \right.$$
$$\left. \times (1 - 2(\gamma+1)l + 4(\gamma+1)^2 l^2) \right\}^{-1},$$

and are able to find a positive number $\varepsilon'' > 0$, such that for $0 < l \leq 2\frac{\gamma}{\gamma+1}$,

$$-8\frac{(\gamma+1)^4}{\gamma^3} l^2$$
$$\times \frac{1 - (\gamma-1)(\gamma+1)/\gamma l + 4((\gamma+1)/\gamma l)^2}{(1 + 2(\gamma+1)/\gamma l + 4((\gamma+1)/\gamma l)^2)(1 - 2(\gamma+1)l + 4(\gamma+1)^2 l^2)}$$
$$< 3\varepsilon'' l^2,$$

and on the two left-most points of $\Gamma^M_{\text{remote}}$, $(1 + \sqrt{3}i)(2M)^{-10/39}$ and $(1 - \sqrt{3}i)(2M)^{-10/39}$,

$$\Re\left( \frac{(1+\gamma)^4}{3\gamma^3} 2Mw^3 + 2MR(w) \right)\bigg|_{w=(1\pm\sqrt{3}i)M^{-10/39}}$$

$$= 2M\left(-\frac{8}{3}\frac{(\gamma+1)^4}{\gamma^3}(2M)^{-10/13} + O(M^{-40/39})\right)$$

$$= -\frac{8}{3}\frac{(\gamma+1)^4}{\gamma^3}(2M)^{3/13}(1 + O(M^{-1/39}))$$

$$< 2M \int_0^{(2M)^{-10/39}} -3\varepsilon'' t^2 \, dt.$$

Therefore we know that for $w \in \Gamma_{\text{remote}}^M$,

$$\Re\left(\frac{(1+\gamma)^4}{3\gamma^3} 2Mw^3 + 2MR_1(w)\right) < 2M \int_0^l -3\varepsilon'' t^2 \, dt = -2M\varepsilon'' l^3,$$

and have the estimation that if $\eta \geq T$, $0 < \varepsilon''' < \varepsilon''$ and $M$ large enough $[l \geq (2M)^{-10/39}]$,

$$|F_{aM}(\xi, w)| < \frac{(1+\gamma)^{4/3}}{\gamma}(2M)^{1/3}$$

$$\times \exp\left\{-(\eta - T)\frac{(1+\gamma)^{4/3}}{\gamma}(2M)^{1/3} l\right.$$

$$\left. - \left(\varepsilon'' l^3 + T\frac{(1+\gamma)^{4/3}}{\gamma}(2M)^{-2/3} l\right) 2M\right\}$$

$$\times \left|\frac{-(\gamma+1)w + 1}{(\gamma+1)/\gamma w + 1} \frac{1}{w + (a\gamma - 1)/((\gamma+1)(a+1))}\right|$$

$$< e^{-(\eta - T)(1+\gamma)^{4/3}/\gamma (2M)^{1/13} - \varepsilon'''(2M)^{3/13}}.$$

Now we get the result by inequalities similar to (53)–(54). □

LEMMA 4. *If $T$ is fixed and $c$ is large enough,*

$$\left|\frac{1}{2\pi i}\int_{\Gamma_{\geq c}^\infty} e^{-Tu + u^3/3} \, du\right| < \frac{1}{c}.$$

PROOF. Obvious. □

LEMMA 5. *If $T$ is fixed and $M$ is large enough, then for any $\eta > T$, holds:*

$$\left|\frac{1}{2\pi i}\int_{\Gamma_{\text{local}}^M} F_{aM}(\eta, w) \, dw - \frac{(\gamma+1)(a+1)}{1 - a\gamma}\operatorname{Ai}(\eta)\right| < \frac{1}{3}\frac{e^{-\eta/2}}{M^{1/40}}.$$

PROOF. On $\Gamma_{\text{local}}^M$, $|w| \leq 2(2M)^{10/39}$, so by (49)

$$F_{aM}(\eta, w) = \frac{(\gamma+1)(a+1)}{a\gamma - 1} \frac{(1+\gamma)^{4/3}}{\gamma} (2M)^{1/3}$$
$$\times e^{-(1+\gamma)^{4/3}/\gamma(2M)^{1/3}\eta w + (1+\gamma)^4/(3\gamma^3)2Mw^3}(1 + O(M^{-1/39})).$$

After the substitution $u = \frac{(1+\gamma)^{4/3}}{\gamma}(2M)^{1/3}w$, we get

$$\int_{\Gamma_{\text{local}}^M} F_{aM}(\eta, w)\, dw = \frac{(\gamma+1)(a+1)}{(a\gamma-1)} \int_{\Gamma_{<(1+\gamma)^{4/3}/\gamma(2M)^{1/13}}^\infty} e^{-\eta u + u^3/3}\, du$$
(55)
$$\times (1 + O(M^{-1/39})),$$

and the $O(M^{-1/39})$ term is independent to $w$.

On $\Gamma^\infty$, if $\eta > T$, $e^{-(\eta-T)u} < e^{T/2}e^{-\eta/2}$. By (2) and (55), we have

$$\left| \frac{1}{2\pi i} \int_{\Gamma_{\text{local}}^M} F_{aM}(\eta, w)\, dw - \frac{(\gamma+1)(a+1)}{1 - a\gamma}\operatorname{Ai}(\eta) \right|$$
$$< e^{T/2}e^{-\eta/2} \left| \frac{(\gamma+1)(a+1)}{2\pi i(1-a\gamma)} \int_{\Gamma_{\geq(1+\gamma)^{4/3}/\gamma(2M)^{1/13}}^\infty} |e^{-Tu+u^3/3}|\, du \right|$$
$$+ e^{T/2}e^{-\eta/2} \left| \frac{(\gamma+1)(a+1)}{2\pi i(1-a\gamma)} \right.$$
$$\left. \times \int_{\Gamma_{<(1+\gamma)^{4/3}/\gamma(2M)^{1/13}}^\infty} |e^{-Tu+u^3/3}|\, du\, O(M^{-1/39}) \right|,$$

and we can get the result by direct calculation. $\square$

CONCLUSION OF THE PROOF OF (43). Putting Lemmas 2–5 together, we get the convergence in $L^2$ norm:

$$\lim_{M\to\infty} (1+a)^{2(N-M)-1}a^{-2N+1}\frac{(\gamma+1)^{2(M-N)+1/3}}{\gamma^{2M}}(1-a\gamma)(2M)^{1/3}$$
$$\times e^{-2M/(1+\gamma)y}\varphi_{2N-1}(y)\chi_T(\eta) = \operatorname{Ai}(\eta)\chi_T(\eta).$$

On the other hand, for $\eta \in [T, \infty)$,

(56) $\quad \lim_{M\to\infty} (1+\gamma^{-1})^{2(N-M)-1}e^{M-N}y^{M-N+1/2}e^{(1-\gamma)/(1+\gamma)My} = 1$

and

(57) $\quad (1+\gamma^{-1})^{2(N-M)-1}e^{M-N}y^{M-N+1/2}e^{-(1-\gamma)/(1+\gamma)My} \leq 1 + O\left(\frac{\eta}{\sqrt{M}}\right).$

Therefore, in $L^2$ norm,

$$\lim_{M\to\infty}(1+a)^{2(N-M)-1}a^{-2N+1}\frac{(1-a\gamma)(2M)^{1/3}}{(\gamma+1)^{2/3}\gamma^{2N-1}}e^{M-N}\psi_{2N-1}(y)\chi(\eta)$$
$$=\operatorname{Ai}(\eta)\chi(\eta). \qquad \square$$

Now we conclude the proof of the $-1<a<\gamma^{-1}$ part of Theorem 1. By Stirling's formula, we get

(58) $$\lim_{M\to\infty}(2M)^{2(M-N)}\frac{(2N-1)!}{(2M-1)!}e^{2(N-M)}\gamma^{4N-1}=1,$$

and then by (34), (38) and Proposition 5, we have the convergence in trace norm

(59)
$$\lim_{M\to\infty}\frac{(1-a\gamma)(2M)^{1/3}}{(1+\gamma)^{2/3}}\chi(\xi)\widetilde{S}_{4b}(\xi,\eta)\chi(\eta)$$
$$=\frac{1}{2}\chi(\xi)\Big(\operatorname{Ai}(\xi)\operatorname{Ai}(\eta)+\operatorname{Ai}'(\xi)\int_\eta^\infty\operatorname{Ai}(t)\,dt\Big)\chi(\eta),$$

which implies that in trace norm,

$$\lim_{M\to\infty}\chi(\xi)\widetilde{S}_{4b}(\xi,\eta)\chi(\eta)=0.$$

Now we get the desired result

$$\lim_{M\to\infty}\chi(\xi)\widetilde{S}_4(\xi,\eta)\chi(\eta)=\lim_{M\to\infty}\chi(\xi)\widetilde{S}_{4a}(\xi,\eta)\chi(\eta)=\chi(\xi)\widehat{S}_4(\xi,\eta)\chi(\eta),$$

and in the same way

$$\lim_{M\to\infty}\chi(\xi)\widetilde{SD}_4(\xi,\eta)\chi(\eta)=\chi(\xi)\widehat{SD}_4(\xi,\eta)\chi(\eta),$$
$$\lim_{M\to\infty}\chi(\xi)\widetilde{IS}_4(\xi,\eta)\chi(\eta)=\chi(\xi)\widehat{IS}_4(\xi,\eta)\chi(\eta).$$

Therefore, in trace norm

$$\lim_{M\to\infty}\widetilde{P}_T(\xi,\eta)\chi(\eta)=\chi(\xi)\begin{pmatrix}\widehat{S}_4(\xi,\eta) & \widehat{SD}_4(\xi,\eta) \\ \widehat{IS}_4(\xi,\eta) & \widehat{S}_4(\eta,\xi)\end{pmatrix}\chi(\eta),$$

and the convergence of Fredholm determinant follows.

3.2. *Proof of the $a=\gamma^{-1}$ part of Theorem 1.* When $a=\gamma^{-1}$, the $1-a\gamma^{-1}$ in (59) vanishes, so we need other asymptotic formulas for $\psi_{2N-1}(\eta)$ and $\psi'_{2N-1}(\eta)$. The approach is similar to that in the $a<\gamma^{-1}$ case, so we just sketch the proof.

PROPOSITION 6.  *For fixed $\gamma \geq 1$, $a = \gamma^{-1}$, $\varepsilon > 0$ and any $T$, we have the convergences in $L^2$ norm with respect to $\xi$ or $\eta$:*

$$\lim_{M \to \infty} \gamma^{-2N-1}(1+\gamma)^{4/3}(2M)^{1/3} e^{M-N} e^{\varepsilon \xi} L_{2N-1}^{(2(M-N))}$$
$$\times (2Mx) x^{M-N-1/2} e^{-Mx} \chi(\xi)$$
$$= -e^{\varepsilon \xi} \operatorname{Ai}(\xi) \chi(\xi),$$

$$\lim_{M \to \infty} \gamma^{-2N} 2M e^{M-N} e^{-\varepsilon \eta} \int_y^\infty L_{2N-1}^{(2(M-N))}(2Mt) t^{M-N-1/2} e^{-Mt}\, dt\, \chi(\eta)$$

(60) $$= -e^{-\varepsilon \eta} \int_\eta^\infty \operatorname{Ai}(t)\, dt\, \chi(\eta),$$

$$\lim_{M \to \infty} (1+a)^{2(N-M)-1} a^{-2N+1} e^{M-N} \gamma^{-2N+1} e^{-\varepsilon \eta} \psi_{2N-1}(y) \chi(\eta)$$
$$= e^{-\varepsilon \eta} s^{(1)}(\eta) \chi(\eta),$$

$$\lim_{M \to \infty} (1+a)^{2(N-M)-1} a^{-2N+1} e^{M-N} \frac{(\gamma+1)^{4/3}}{\gamma^{2N}(2M)^{-2/3}} e^{\varepsilon \xi} \psi'_{2N-1}(x) \chi(\xi)$$
$$= e^{\varepsilon \xi} \operatorname{Ai}(\xi) \chi(\xi).$$

SKETCH OF PROOF OF (60).  We perform the same algebraic procedure and use the contour $\bar{\bar{\Gamma}}^M = \bar{\bar{\Gamma}}_1^M \cup \bar{\bar{\Gamma}}_2^M \cup \bar{\bar{\Gamma}}_3^M \cup \bar{\bar{\Gamma}}_4^M$ which is slightly different from the $\Gamma^M$ in the $a < \gamma^{-1}$ case (see Figure 4):

$$\bar{\bar{\Gamma}}_1^M = \left\{ (4-t) \frac{\gamma}{\gamma+1} e^{\pi i/3} \Big| 0 \leq t \leq 4 - \frac{\varepsilon/2}{(1+\gamma)^{1/3}} (2M)^{-1/3} \right\},$$

$$\bar{\bar{\Gamma}}_2^M = \left\{ \frac{\gamma \varepsilon/2}{(1+\gamma)^{4/3}} (2M)^{-1/3} e^{t \pi i} \Big| \frac{1}{3} \leq t \leq \frac{5}{3} \right\},$$

$$\bar{\bar{\Gamma}}_3^M = \left\{ t \frac{\gamma}{\gamma+1} e^{5\pi i/3} \Big| \frac{\varepsilon/2}{(1+\gamma)^{1/3}} (2M)^{-1/3} \leq t \leq 4 \right\},$$

$$\bar{\bar{\Gamma}}_4^M = \left\{ 2 \frac{\gamma}{\gamma+1} + it \Big| -2\sqrt{3} \frac{\gamma}{\gamma+1} \leq t \leq 2\sqrt{3} \frac{\gamma}{\gamma+1} \right\}$$

and for asymptotic analysis, we define $\bar{\bar{\Gamma}}_{\text{remote}}^M$, $\bar{\bar{\Gamma}}_{\text{local}}^M$, $\bar{\bar{\Gamma}}_{<c}^\infty$ and $\bar{\bar{\Gamma}}_{\geq c}^\infty$ in the same way as (50)–(51). Then we get

$$\varphi_{2N-1}(y) = -(1+a)^{2(M-N)+1} a^{2N-1}$$
$$\times \frac{\gamma^{2M}}{(\gamma+1)^{2(M-N)+1}} e^{(2M)/(1+\gamma)y}$$
$$\times \frac{1}{2\pi i} \oint_{\bar{\bar{\Gamma}}^M} \exp\left\{ -\frac{(1+\gamma)^{4/3}}{\gamma}(2M)^{1/3} \eta w \right.$$

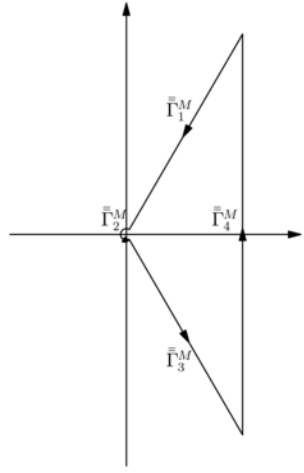

Fig. 4. $\bar{\bar{\Gamma}}^M$.

$$+ \frac{(1+\gamma)^4}{3\gamma^3} 2Mw^3 + 2MR_1(w)\bigg\}$$
$$\times \frac{-(\gamma+1)w+1}{(\gamma+1)/\gamma w+1} \frac{dw}{w}.$$

If we denote

$$F_M(\eta, w) = \exp\bigg\{-\frac{(1+\gamma)^{4/3}}{\gamma}(2M)^{1/3}\eta w + \frac{(1+\gamma)^4}{3\gamma^3} 2Mw^3 + 2MR_1(w)\bigg\}$$
$$\times \frac{-(\gamma+1)w+1}{(\gamma+1)/\gamma w+1} \frac{1}{w},$$

then parallel to Lemmas 2–5, we have:

LEMMA 6. *For any $T$ fixed, and $M$ large enough, if $\eta > T$, then*
$$\bigg|e^{-\varepsilon\eta} \frac{1}{2\pi i} \int_{\bar{\bar{\Gamma}}_4^M} F_M(\eta, w)\, dw\bigg| < \frac{1}{3} \frac{e^{-\varepsilon\eta/2}}{M^{1/40}}.$$

LEMMA 7. *For any $T$ fixed, and $M$ large enough, if $\eta > T$, then*
$$\bigg|e^{-\varepsilon\eta} \frac{1}{2\pi i} \int_{\bar{\bar{\Gamma}}_{\text{remote}}^M} F_M(\eta, w)\, dw\bigg| < \frac{1}{3} \frac{e^{-\varepsilon\eta/2}}{M^{1/40}}.$$

LEMMA 8. *If $T$ is fixed and $c$ is large enough,*
$$\bigg|\frac{1}{2\pi i} \int_{\bar{\bar{\Gamma}}_{\geq c}^\infty} e^{-Tu + u^3/3} \frac{du}{u}\bigg| < \frac{1}{c}.$$

LEMMA 9. *For any $T$ fixed, and $M$ large enough, if $\eta > T$, then*

$$\left| e^{-\varepsilon\eta} \frac{1}{2\pi i} \int_{\bar{\Gamma}^M_{\text{local}}} F_M(\eta, w) \, dw - e^{-\eta/2} s^{(1)}(\eta) \right| < \frac{1}{3} \frac{e^{-\varepsilon\eta/2}}{M^{1/40}}.$$

Using Lemmas 6–9, we get the convergence in $L^2$ norm:

$$\lim_{M \to \infty} (1+a)^{2(N-M)-1} a^{-2N+1} \frac{(\gamma+1)^{2(M-N)+1}}{\gamma^{2M}} e^{-2M/(1+\gamma)y}$$
$$\times e^{-\varepsilon\eta} \varphi_{2N-1}(y) \chi(\eta)$$
$$= e^{-\varepsilon\eta} s^{(1)}(\eta) \chi(\eta).$$

Furthermore, because of the limit result (56) and (57), we get the $L^2$ convergence

$$\lim_{M \to \infty} (1+a)^{2(N-M)-1} a^{-2N+1} \gamma^{-2N+1} e^{M-N} e^{-\varepsilon\eta} \psi_{2N-1}(y) \chi(\eta)$$
$$= e^{-\varepsilon\eta} s^{(1)}(\eta) \chi(\eta). \qquad \square$$

Now we conclude the proof of the $a > \gamma^{-1}$ part of Theorem 1. Using (34), (58) and Proposition 6 we have the convergence in trace norm

$$\lim_{M \to \infty} \chi(\xi) e^{\varepsilon\xi} \widetilde{S}_{4b}(\xi, \eta) e^{-\varepsilon\eta} \chi(\eta)$$
$$= \frac{1}{2} \chi(\xi) e^{\varepsilon\xi} \left( \text{Ai}(\xi) s^{(1)}(\eta) + \text{Ai}(\xi) \int_\eta^\infty \text{Ai}(t) \, dt \right) e^{-\varepsilon\eta} \chi(\eta)$$
$$= \frac{1}{2} \chi(\xi) e^{\varepsilon\xi} \text{Ai}(\xi) e^{-\varepsilon\eta} \chi(\eta),$$

and this together with the conjugated convergence result (discussed in the Appendix) of $\widetilde{S}_{4a}(\xi, \eta)$ in formulas (39) and (40) of Section 3.1 conclude

(61) $$\lim_{M \to \infty} \chi(\xi) e^{\varepsilon\xi} \widetilde{S}_4(\xi, \eta) e^{-\varepsilon\eta} \chi(\eta) = \chi(\xi) e^{\varepsilon\xi} \overline{\overline{S}}_4(\xi, \eta) e^{-\varepsilon\eta} \chi(\eta).$$

In the same way we get

$$\lim_{M \to \infty} \chi(\xi) e^{\varepsilon\xi} \widetilde{SD}_4(\xi, \eta) e^{\varepsilon\eta} \chi(\eta) = \chi(\xi) e^{\varepsilon\xi} \overline{\overline{SD}}_4(\xi, \eta) e^{\varepsilon\eta} \chi(\eta),$$

$$\lim_{M \to \infty} \chi(\xi) e^{-\varepsilon\xi} \widetilde{IS}_4(\xi, \eta) e^{-\varepsilon\eta} \chi(\eta) = \chi(\xi) e^{-\varepsilon\xi} \overline{\overline{IS}}_4(\xi, \eta) e^{-\varepsilon\eta} \chi(\eta).$$

Then we get the convergence in trace norm of a conjugate of $\widetilde{P}_T(\xi, \eta)$

$$\lim_{M \to \infty} \chi(\xi) \begin{pmatrix} e^{\varepsilon\xi} \widetilde{S}_4(\xi, \eta) e^{-\varepsilon\eta} & e^{\varepsilon\xi} \widetilde{SD}_4(\xi, \eta) e^{\varepsilon\eta} \\ e^{-\varepsilon\xi} \widetilde{IS}_4(\xi, \eta) e^{-\varepsilon\eta} & e^{-\varepsilon\xi} \widetilde{S}_4(\eta, \xi) e^{\varepsilon\eta} \end{pmatrix} \chi(\eta)$$
$$= \chi(\xi) \begin{pmatrix} e^{\varepsilon\xi} \overline{\overline{S}}_4(\xi, \eta) e^{-\varepsilon\eta} & e^{\varepsilon\xi} \overline{\overline{SD}}_4(\xi, \eta) e^{\varepsilon\eta} \\ e^{-\varepsilon\xi} \overline{\overline{IS}}_4(\xi, \eta) e^{-\varepsilon\eta} & e^{-\varepsilon\xi} \overline{\overline{S}}_4(\eta, \xi) e^{\varepsilon\eta} \end{pmatrix} \chi(\eta),$$

and the convergence of Fredholm determinant follows.

3.3. *Proof of the $a > \gamma^{-1}$ part of Theorem 1.* If $a > \gamma^{-1}$, the location as well as the fluctuation scale of the largest sample eigenvalue is changed. We change variables as $p_M = (a+1)(1 + \frac{1}{\gamma^2 a})$ and $q_M = (a+1)\sqrt{1 - \frac{1}{\gamma^2 a^2}} \frac{1}{\sqrt{2M}}$, and then by (36) the kernel $S_*(x,y)$ after substitution is [here $*$ stands for 4, 4a or 4b, and the $\widetilde{S}_*(\xi, \eta)$ in this subsection is not identical to that in Sections 3.1 and 3.2]

$$\widetilde{S}_*(\xi, \eta) = (a+1)\sqrt{1 - \frac{1}{\gamma^2 a^2}} \frac{1}{\sqrt{2M}}$$
(62)
$$\times S_*(x,y)\Big|_{\substack{x=(a+1)(1+1/(\gamma^2 a))+(a+1)\sqrt{1-1/(\gamma^2 a^2)}1/(\sqrt{2M})\xi \\ y=(a+1)(1+1/(\gamma^2 a))+(a+1)\sqrt{1-1/(\gamma^2 a^2)}1/(\sqrt{2M})\eta}}.$$

We analyze $\widetilde{S}_{4b}(\xi, \eta)$ first.

PROPOSITION 7. *For fixed $\gamma \geq 1$, $a > \gamma^{-1}$, $\varepsilon > 0$ and any $T$, we have convergences in $L^2$ norm with respect to $\xi$ or $\eta$:*

$$\lim_{M \to \infty} \frac{(\gamma^2 a + 1)^{M-N+1/2}}{(\gamma^2 a)^{M+N+1/2}(a+1)^{M-N-1/2}} \sqrt{(\gamma^2 a^2 - 1)2M} e^{M-N}$$

$$\times e^{(\gamma^2 a^2 - 1)/((\gamma^2 a + 1)(a+1))Mx}$$
(63)
$$\times e^{\varepsilon \xi} L_{2N-1}^{(2(M-N))}(2Mx) x^{M-N-1/2} e^{-Mx} \chi(\xi)$$

$$= -\frac{1}{\sqrt{2\pi}} \exp\left\{-\frac{1}{4} \frac{\gamma^4 a^2 + \gamma^2 a^2 + 4\gamma^2 a + \gamma^2 + 1}{(\gamma^2 a + 1)^2} \xi^2 + \varepsilon \xi\right\} \chi(\xi),$$

$$\lim_{M \to \infty} \frac{1}{2} \frac{(\gamma^2 a + 1)^{M-N-1/2}(\gamma^2 a^2 - 1)}{(\gamma^2 a)^{M+N+1/2}(a+1)^{M-N+1/2}} \sqrt{\gamma^2 a^2 - 1}(2M)^{3/2} e^{M-N}$$

$$\times e^{(\gamma^2 a^2 - 1)/((\gamma^2 a + 1)(a+1))My} e^{\varepsilon \eta}$$
(64)
$$\times \int_y^\infty L_{2N-1}^{(2(M-N))}(2Mt) t^{M-N-1/2} e^{-Mt} dt \, \chi(\eta)$$

$$= -\frac{1}{\sqrt{2\pi}} e^{-1/4(\gamma^4 a^2 + \gamma^2 a^2 + 4\gamma^2 a + \gamma^2 + 1)/(\gamma^2 a + 1)^2 \eta^2 + \varepsilon \eta} \chi(\eta),$$

$$\lim_{M \to \infty} \left(\frac{\gamma^2 a}{(\gamma^2 a + 1)(a+1)}\right)^{M-N+1/2} e^{M-N}$$

(65)
$$\times e^{-(\gamma^2 a^2-1)/((\gamma^2 a+1)(a+1))My} e^{-\varepsilon \eta} \psi_{2N-1}(y)\chi(\eta)$$
$$= e^{-1/4(\gamma^2 a^2-1)(\gamma^2-1)/(\gamma^2 a+1)^2 \eta^2 - \varepsilon \eta} \chi(\eta),$$

(66)
$$\lim_{M\to\infty} \left(\frac{\gamma^2 a}{(\gamma^2 a+1)(a+1)}\right)^{M-N-1/2} e^{M-N} \frac{\gamma^2 a}{(\gamma^2 a^2-1)M}$$
$$\times e^{-(\gamma^2 a^2-1)/((\gamma^2 a+1)(a+1))Mx} e^{-\varepsilon \xi} \psi'_{2N-1}(x)\chi(\xi)$$
$$= e^{-1/4(\gamma^2 a^2-1)(\gamma^2-1)/(\gamma^2 a+1)^2 \xi^2 - \varepsilon \xi} \chi(\xi).$$

We only prove (65). The identity (45) still holds, but we need to use another contour and a new procedure of steepest-descent analysis.

Since
$$e^{-2Myz} \frac{(z+1)^{2M}}{z^{2N}}$$
$$= e^{2M(-(a+1)(1+1/(\gamma^2 a))z + \log(z+1) - \gamma^{-2}\log z) - (a+1)\sqrt{1-1/(\gamma^2 a^2)}\sqrt{2M}\eta z},$$

if we denote (ignoring the ambiguity of values of logarithm)
$$g(z) = -(a+1)\left(1 + \frac{1}{\gamma^2 a}\right) z + \log(z+1) - \gamma^{-2} \log z,$$

then we get:

- $g'(z) = -(a+1)(1 + \frac{1}{\gamma^2 a}) + \frac{1}{z+1} - \frac{\gamma^{-2}}{z}$, with zero points $z = -\frac{1}{1+\gamma^2 a}$ and $z = -\frac{a}{1+a}$;
- $g''(z) = -\frac{1}{(z+1)^2} + \frac{\gamma^{-2}}{z^2}$, $g''(-\frac{1}{1+\gamma^2 a}) = (\gamma^{-1}+\gamma a)^2(1-\frac{1}{\gamma^2 a^2}) > 0$ and $g''(-\frac{a}{1+a}) = (1+a)^2(\frac{1}{\gamma^2 a^2}-1) < 0$.

So we take $z = -\frac{1}{1+\gamma^2 a}$ as the saddle point, and locally around that point, after the substitution $w = z + \frac{1}{1+\gamma^2 a}$, we get
$$g\left(-\frac{1}{1+\gamma^2 a} + w\right) = \frac{a+1}{\gamma^2 a} + \log(\gamma^2 a) - (1-\gamma^{-2})\log(\gamma^2 a+1) + \gamma^{-2}\pi i$$
$$+ \frac{1}{2}(\gamma^{-1}+\gamma a)^2\left(1 - \frac{1}{\gamma^2 a^2}\right)w^2 + R_2(w),$$

where
$$R_2(w) = O(w^3) \quad \text{as } w \to 0,$$

so that
$$\oint_C e^{2M(-yz+\log(z+1)-\gamma^{-2}\log z)} \frac{z}{((a+1)z+a)(z+1)} dz$$

$$\begin{aligned}(67)\quad &= \oint_{\Sigma^M} \exp\Big\{2M\Big(\frac{a+1}{\gamma^2 a} + \log(\gamma^2 a) - (1-\gamma^{-2})\log(\gamma^2 a + 1) \\
&\quad + \gamma^{-2}\pi i + \frac{1}{2}(\gamma^{-1}+\gamma a)^2\Big(1-\frac{1}{\gamma^2 a^2}\Big)w^2 + R_2(w) \\
&\quad - (a+1)\sqrt{1-\frac{1}{\gamma^2 a^2}}\frac{\eta}{\sqrt{2M}}\Big(w-\frac{1}{\gamma^2 a+1}\Big)\Big)\Big\} \\
&\quad \times \frac{w - 1/(\gamma^2 a + 1)}{((a+1)w + (\gamma^2 a^2 - 1)/(\gamma^2 a + 1))(w + (\gamma^2 a)/(\gamma^2 a + 1))}\,dw \\
&= -\frac{1}{a+1}\frac{(\gamma^2 a)^{2M-1}}{(\gamma^2 a + 1)^{2(M-N)}}e^{2M/(\gamma^2 a+1)x} \\
&\quad \times \oint_{\Sigma^M} \exp\Big\{-(a+1)\sqrt{1-\frac{1}{\gamma^2 a^2}}\sqrt{2M}\eta w \\
&\quad + \frac{1}{2}(\gamma^{-1}+\gamma a)^2\Big(1-\frac{1}{\gamma^2 a^2}\Big)2Mw^2 + 2MR_2(w)\Big\} \\
&\quad \times \frac{-(\gamma^2 a + 1)w + 1}{(\gamma^2 a + 1)/(\gamma^2 a)w + 1} \\
&\quad \times \frac{1}{w + (\gamma^2 a^2 - 1)/((\gamma^2 a + 1)(a+1))}\,dw,\end{aligned}$$

where $\Sigma^M$ is a contour around $\frac{1}{\gamma^2 a+1}$, composed of $\Sigma_1^M$, $\Sigma_2^M$, $\Sigma_3^M$ and $\Sigma_4^M$, which are defined as (see Figure 5)

$$\Sigma_1^M = \{-it|-2 \leq t \leq 2\}, \qquad \Sigma_2^M = \{4-t+2i|0 \leq t \leq 4\},$$
$$\Sigma_3^M = \{4+it|-2 \leq t \leq 2\}, \qquad \Sigma_4^M = \{t-2i|0 \leq t \leq 4\}.$$

And for the asymptotic analysis, we define (see Figure 6)

$$\Sigma_{\text{local}}^M = \{w \in \Sigma^M||w| \leq M^{-2/5}\}, \qquad \Sigma_{\text{remote}}^M = \Sigma_1^M \setminus \Sigma_{\text{local}}^M,$$

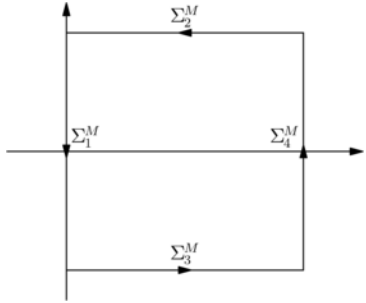

FIG. 5. $\Sigma^M$.

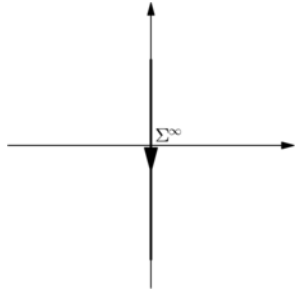

Fig. 6. $\Sigma^{\infty}$.

$$\Sigma^{\infty} = \{-it | -\infty < t < \infty\}, \qquad \Sigma^{\infty}_{<c} = \{w \in \Sigma^{\infty} | |w| < c\},$$
$$\Sigma^{\infty}_{\geq c} = \Sigma^{\infty} \setminus \Sigma^{\infty}_{<c}.$$

Then if we denote

$$G_{aM}(\eta, w) = (\gamma^{-1} + \gamma a)\sqrt{\left(1 - \frac{1}{\gamma^2 a^2}\right)2M}$$

$$\times \exp\Bigg\{-(a+1)\sqrt{\left(1 - \frac{1}{\gamma^2 a^2}\right)2M}\eta w$$

$$+ \frac{1}{2}(\gamma^{-1} + \gamma a)^2 \left(1 - \frac{1}{\gamma^2 a^2}\right) 2Mw^2 + 2MR_2(w)\Bigg\}$$

$$\times \frac{-(\gamma^2 a + 1)w + 1}{(\gamma^2 a + 1)/(\gamma^2 a)w + 1} \frac{1}{w + (\gamma^2 a^2 - 1)/((\gamma^2 a + 1)(a+1))},$$

we have four lemmas similar to Lemmas 2–5:

LEMMA 10. *For any $T$ fixed, and $M$ large enough, if $\eta > T$, then*

$$\left| e^{-\varepsilon\eta} \frac{1}{2\pi i} \int_{\Sigma^M_2 \cup \Sigma^M_3 \cup \Sigma^M_4} G_{aM}(\eta, w)\, dw \right| < \frac{1}{3} \frac{e^{-\varepsilon\eta}}{M^{1/10}}.$$

LEMMA 11. *For any $T$ fixed, and $M$ large enough, if $\eta > T$, then*

$$\left| e^{-\varepsilon\eta} \frac{1}{2\pi i} \int_{\Sigma^M_{\text{remote}}} G_{aM}(\eta, w)\, dw \right| < \frac{1}{3} \frac{e^{-\varepsilon\eta}}{M^{1/10}}.$$

LEMMA 12. *If $T$ is fixed and $c$ is large enough,*

$$\left| \frac{1}{2\pi i} \int_{\Gamma^{\infty}_{\geq c}} e^{-(a+1)/(\gamma^{-1} + \gamma a)Tu + u^2/2}\, du \right| < \frac{1}{c}.$$

LEMMA 13. *For any $T$ fixed, and $M$ large enough, if $\eta > T$, then*

$$\left| e^{-\varepsilon\eta} \frac{1}{2\pi i} \int_{\Sigma_{\text{local}}^M} G_{aM}(\eta, w)\, dw \right.$$
$$\left. - \frac{(\gamma^2 a + 1)(a+1)}{(\gamma^2 a^2 - 1)\sqrt{2\pi}} e^{-1/2((\gamma(a+1))/(\gamma^2 a+1)\eta)^2 - \varepsilon\eta} \right|$$
$$< \frac{1}{3} \frac{e^{-\varepsilon\eta}}{M^{1/10}}.$$

Their proofs are the same as those of Lemmas 2–5, and we need the identity

$$\frac{1}{2\pi i} \int_{\Sigma^\infty} e^{-(a+1)/(\gamma^{-1}+\gamma a)\eta u + u^2/2}\, du = -\frac{1}{\sqrt{2\pi}} e^{-1/2(\gamma(a+1)/(\gamma^2 a+1)\eta)^2}.$$

SKETCH OF PROOF OF (65). Because the pole $z = -\frac{a}{a+1}$, which is $w = -\frac{\gamma^2 a^2 - 1}{(\gamma^2 a + 1)(a+1)}$ in the $w$ plane, is not in side of $\Sigma^M$, so

(68) $$\oint_C e^{-2Mxz} \frac{(z+1)^{2(M-N)}}{z + a/(a+1)}\, dz = 0.$$

Similar to but subtler than (56) and (57), if we denote (here we have a notation conflict with the $r_i$ defined in Section 2.3, but there should be no confusion)

$$r_M(\eta) = \left( \frac{\gamma^2 a}{(\gamma^2 a + 1)(a+1)} \right)^{M-N+1/2} e^{M-N} y^{M-N+1/2}$$
$$\times e^{-(\gamma^2-1)a/((\gamma^2 a+1)(a+1))My - \varepsilon\eta},$$

we have for $\eta \in [T, \infty)$,

$$\lim_{M \to \infty} r_M(\eta) = e^{-1/4(\gamma^2 a^2 - 1)(\gamma^2 - 1)/(\gamma^2 a + 1)^2 \eta^2 - \varepsilon\eta},$$

and for a large enough positive $C$, $\eta \in [C, \infty)$ and $M_1 < M_2$, pointwisely

$$r_{M_1}(\eta) > r_{M_2}(\eta) > 0,$$

so that we can use the dominant convergence theorem to prove that in $L^2$ norm,

$$\lim_{M \to \infty} r_M(\eta)\chi(\eta) = e^{-1/4(\gamma^2 a^2 - 1)(\gamma^2 - 1)/((\gamma^2 a + 1)^2)\eta^2 - \varepsilon\eta}\chi(\eta).$$

Finally, since from (45), (67) and (68),

$$\psi_{2N-1}(y) = y^{M-N+1/2} e^{(a-1)/(a+1)My}$$

$$+ (1+a)^{2(M-N)} a^{2N-1}$$
$$\times \frac{(\gamma^2 a)^{2M}}{(\gamma^2 a + 1)^{2(M-N)+1}} \frac{1}{\sqrt{(\gamma^2 a^2 - 1)2M}}$$
$$\times y^{M-N+1/2} e^{(1-\gamma^2 a)/(1+\gamma^2 a)My} \frac{1}{2\pi i} \int_{\Sigma^M} G_{aM}(\eta, w) \, dw$$
$$= y^{M-N+1/2} e^{-((\gamma^2-1)a)/((\gamma^2 a+1)(a+1))My}$$
$$\times \Bigg[ e^{(\gamma^2 a^2 - 1)/(\gamma^2 a + 1)(a+1)My}$$
$$+ (1+a)^{2(M-N)} a^{2N-1} \frac{(\gamma^2 a)^{2M}}{(\gamma^2 a + 1)^{2(M-N)+1}} \frac{1}{\sqrt{(\gamma^2 a^2 - 1)2M}}$$
$$\times e^{-(\gamma^2 a^2 - 1)/((\gamma^2 a + 1)(a+1))My} \frac{1}{2\pi i} \int_{\Sigma^M} G_{aM}(\eta, w) \, dw \Bigg],$$

we get

$$\left( \frac{\gamma^2 a}{(\gamma^2 a + 1)(a+1)} \right)^{M-N+1/2}$$
$$\times e^{M-N} e^{-(\gamma^2 a^2 - 1)/((\gamma^2 a + 1)(a+1))My} e^{-\varepsilon \eta} \psi_{2N-1}(y) \chi(\eta)$$
$$= r_M(\eta) \Bigg[ 1 + (1+a)^{2(M-N)} a^{2N-1} \frac{(\gamma^2 a)^{2M}}{(\gamma^2 a + 1)^{2(M-N)+1}}$$
$$\times \frac{1}{\sqrt{(\gamma^2 a^2 - 1)2M}} \exp\left\{ -\frac{(\gamma^2 a^2 - 1)}{(\gamma^2 a + 1)(a+1)} 2My \right\}$$
$$\times \frac{1}{2\pi i} \int_{\Sigma^M} G_{aM}(\eta, w) \, dw \Bigg] \chi(\eta),$$

and get the $L^2$ convergence

$$\lim_{M \to \infty} \left( \frac{\gamma^2 a}{(\gamma^2 a + 1)(a+1)} \right)^{M-N+1/2}$$
$$\times e^{M-N} e^{-(\gamma^2 a^2 - 1)/((\gamma^2 a + 1)(a+1))My} e^{-\varepsilon \eta} \psi_{2N-1}(y) \chi(\eta)$$
$$= \lim_{M \to \infty} r_M(\eta) \chi(\eta) = e^{-1/4(\gamma^2 a^2 - 1)(\gamma^2 - 1)/(\gamma^2 a + 1)^2 \eta^2 - \varepsilon \eta} \chi(\eta),$$

because for $a > \gamma^{-1}$ and $\eta \in [T, \infty)$, we can verify by by elementary but tricky estimation that

$$\lim_{M \to \infty} (1+a)^{2(M-N)} a^{2N-1}$$
$$\times \frac{(\gamma^2 a)^{2M}}{(\gamma^2 a + 1)^{2(M-N)+1}} \frac{e^{-(\gamma^2 a^2 - 1)/((\gamma^2 a + 1)(a+1))2My}}{\sqrt{(\gamma^2 a^2 - 1)2M}} = 0$$

uniformly, and by Lemmas 10–13, in $L^2$ norm

$$\lim_{M\to\infty} e^{-\varepsilon\eta} \frac{1}{2\pi i} \int_{\Sigma^M} G_{aM}(\eta, w)\, dw\, \chi(\eta)$$
$$= \frac{(\gamma^2 a + 1)(a+1)}{(\gamma^2 a^2 - 1)\sqrt{2\pi}} e^{-1/2(\gamma(a+1)/(\gamma^2 a+1)\eta)^2 - \varepsilon\eta} \chi(\eta). \qquad \square$$

For notational simplicity, we denote functions on the left-hand sides of (63)–(66) by $F_1(\xi)\chi(\xi)$, $F_2(\eta)\chi(\eta)$, $F_3(\eta)\chi(\eta)$ and $F_4(\xi)\chi(\xi)$, and denote

$$c_M = (2M)^{2(M-N)} \frac{(2N-1)!}{(2M-1)!} e^{2(N-M)} \gamma^{4N-1}.$$

By (58), we have

$$\lim_{M\to\infty} c_M = 1.$$

Then we get from (34), (62) and (63)–(66)

$$\widetilde{S}_{4b}(\xi, \eta) = -\frac{c_M}{2}(e^{-(\gamma^2 a^2 - 1)/((\gamma^2 a+1)(a+1))M(x-y) - \varepsilon(\xi-\eta)} F_1(\xi) F_3(\eta)$$
(69)
$$+ e^{(\gamma^2 a^2 - 1)/((\gamma^2 a+1)(a+1))M(x-y) + \varepsilon(\xi-\eta)} F_4(\xi) F_2(\eta)).$$

If we define

$$SD_{4a}(x, y) = \sum_{j=0}^{N-2} \frac{1}{r_j}(\psi'_{2j}(x)\psi'_{2j+1}(y) - \psi'_{2j+1}(x)\psi'_{2j}(y)),$$

$$IS_{4a}(x, y) = \sum_{j=0}^{N-2} \frac{1}{r_j}(-\psi_{2j}(x)\psi_{2j+1}(y) + \psi_{2j+1}(x)\psi_{2j}(y)),$$

and

$$SD_{4b}(x, y) = \frac{1}{r_{N-1}}(\psi'_{2N-2}(x)\psi'_{2N-1}(y) - \psi'_{2N-1}(x)\psi'_{2N-2}(y)),$$

$$IS_{4b}(x, y) = \frac{1}{r_{N-1}}(-\psi_{2N-2}(x)\psi_{2N-1}(y) + \psi_{2N-1}(x)\psi_{2N-2}(y)),$$

like

$$S_{4a}(x, y) = \sum_{j=0}^{N-2} \frac{1}{r_j}(-\psi'_{2j}(x)\psi_{2j+1}(y) + \psi'_{2j+1}(x)\psi_{2j}(y)),$$

$$S_{4b}(x, y) = \frac{1}{r_{N-1}}(-\psi'_{2N-2}(x)\psi_{2N-1}(y) + \psi'_{2N-1}(x)\psi_{2N-2}(y)),$$

in (28) and (29), and by (35) and (37) define [$*$ stands for 4, 4a or 4b]

$$\widetilde{SD}_*(\xi,\eta) = (a+1)^2\left(1-\frac{1}{\gamma^2 a^2}\right)\frac{1}{2M}$$
$$\times SD_*(x,y)|_{\substack{x=(a+1)(1+1/(\gamma^2 a))+(a+1)\sqrt{1-1/(\gamma^2 a^2)}1/(\sqrt{2M})\xi,\\ y=(a+1)(1+1/(\gamma^2 a))+(a+1)\sqrt{1-1/(\gamma^2 a^2)}1/(\sqrt{2M})\eta}}$$

$$\widetilde{IS}_*(\xi,\eta) = IS_*(x,y)|_{\substack{x=(a+1)(1+1/(\gamma^2 a))+(a+1)\sqrt{1-1/(\gamma^2 a^2)}1/(\sqrt{2M})\xi,\\ y=(a+1)(1+1/(\gamma^2 a))+(a+1)\sqrt{1-1/(\gamma^2 a^2)}1/(\sqrt{2M})\eta}}$$

like

$$\widetilde{S}_*(\xi,\eta) = (a+1)\sqrt{1-1/(\gamma^2 a^2)}\frac{1}{\sqrt{2M}}$$
$$\times S_*(x,y)|_{\substack{x=(a+1)(1+1/(\gamma^2 a))+(a+1)\sqrt{1-1/(\gamma^2 a^2)}1/(\sqrt{2M})\xi\\ y=(a+1)(1+1/(\gamma^2 a))+(a+1)\sqrt{1-1/(\gamma^2 a^2)}1/(\sqrt{2M})\eta}}$$

in (62), then in the same way of (69), we have

$$\widetilde{SD}_{4b}(\xi,\eta) = \frac{c_M}{4}C_M(e^{-(\gamma^2 a^2-1)/((\gamma^2 a+1)(a+1))M(x-y)-\varepsilon(\xi-\eta)}F_1(\xi)F_4(\eta)$$
$$- e^{(\gamma^2 a^2-1)/((\gamma^2 a+1)(a+1))M(x-y)+\varepsilon(\xi-\eta)}F_4(\xi)F_1(\eta)),$$

$$\widetilde{IS}_{4b}(\xi,\eta) = \frac{c_M}{C_M}(e^{-(\gamma^2 a^2-1)/((\gamma^2 a+1)(a+1))M(x-y)-\varepsilon(\xi-\eta)}F_2(\xi)F_3(\eta)$$
$$- e^{(\gamma^2 a^2-1)/((\gamma^2 a+1)(a+1))M(x-y)+\varepsilon(\xi-\eta)}F_3(\xi)F_2(\eta)),$$

with

$$C_M = \frac{(\gamma^2 a^2-1)^{3/2}\sqrt{2M}}{a\gamma(\gamma^2 a+1)}.$$

Now we write $\widetilde{P}_T(\xi,\eta)$ as the sum

$$\widetilde{P}_T(\xi,\eta) = \widetilde{P}_{Ta}(\xi,\eta) + \widetilde{P}_{Tb}(\xi,\eta),$$

with

$$\widetilde{P}_{Ta}(\xi,\eta) = \chi(\xi)\begin{pmatrix} \widetilde{S}_{4a}(\xi,\eta) & \widetilde{SD}_{4a}(\xi,\eta) \\ \widetilde{IS}_{4a}(\xi,\eta) & \widetilde{S}_{4a}(\eta,\xi) \end{pmatrix}\chi(\eta),$$

$$\widetilde{P}_{Tb}(\xi,\eta) = \chi(\xi)\begin{pmatrix} \widetilde{S}_{4b}(\xi,\eta) & \widetilde{SD}_{4b}(\xi,\eta) \\ \widetilde{IS}_{4b}(\xi,\eta) & \widetilde{S}_{4b}(\eta,\xi) \end{pmatrix}\chi(\eta).$$

If we denote

$$U(\xi) = \begin{pmatrix} \exp\left\{\dfrac{\gamma^2 a^2 - 1}{(\gamma^2 a + 1)(a+1)} M(x-x_0) + \varepsilon\xi\right\} & -\dfrac{C_M}{2}\dfrac{F_4(\xi)}{F_3(\xi)} e^{(\gamma^2 a^2-1)/((\gamma^2 a+1)(a+1))M(x-x_0)+\varepsilon\xi} \\ 0 & e^{(\gamma^2 a^2-1)/((\gamma^2 a+1)(a+1))M(x-x_0)+\varepsilon\xi} \end{pmatrix},$$

$$U^{-1}(\eta) = \begin{pmatrix} e^{-(\gamma^2 a^2-1)/((\gamma^2 a+1)(a+1))M(y-y_0)-\varepsilon\eta} & \dfrac{C_M}{2}\dfrac{F_4(\xi)}{F_3(\xi)} e^{-(\gamma^2 a^2-1)/((\gamma^2 a+1)(a+1))M(y-y_0)-\varepsilon\xi} \\ 0 & e^{(\gamma^2 a^2-1)/((\gamma^2 a+1)(a+1))M(y-y_0)+\varepsilon\xi} \end{pmatrix},$$

with

$$x_0 = y_0 = (a+1)\left(1 + \frac{1}{\gamma^2 a}\right) + (a+1)\sqrt{1 - \frac{1}{\gamma^2 a^2}}\frac{T}{\sqrt{2M}},$$

then we have the result of kernel conjugation

$$U(\xi)\widetilde{P}_{Tb}(\xi,\eta)U^{-1}(\eta) = \chi(\xi)\begin{pmatrix} -\dfrac{c_M}{2}\left(F_1(\xi) + \dfrac{F_2(\xi)F_4(\xi)}{F_3(\xi)}\right)F_3(\eta) & 0 \\ U(\xi)\widetilde{P}_{Tb}(\xi,\eta)U^{-1}(\eta)_{21} & -\dfrac{C_M}{2}F_3(\xi)\left(F_1(\eta) + \dfrac{F_2(\eta)F_4(\eta)}{F_3(\eta)}\right) \end{pmatrix}\chi(\eta),$$

with the entry

$$U(\xi)\widetilde{P}_{Tb}(\xi,\eta)U^{-1}(\eta)_{21}$$
$$= \frac{c_M}{C_M}(e^{-2(\gamma^2 a^2-1)/((\gamma^2 a+1)(a+1))M(x-x_0)-2\varepsilon\xi}F_2(\xi)F_3(\eta)$$
$$- F_3(\xi)F_2(\eta)e^{-2(\gamma^2 a^2-1)/((\gamma^2 a+1)(a+1))M(y-y_0)-2\varepsilon\eta}).$$

We want $U(\xi)\widetilde{P}_{Tb}(\xi,\eta)U^{-1}(\eta)$ to converge in trace norm as $M \to \infty$, and need the results:

LEMMA 14. *In trace norm,*
$$\lim_{M\to\infty} U(\xi)\widetilde{P}_{Tb}(\xi,\eta)U^{-1}(\eta)_{21} = 0.$$

LEMMA 15. *In $L^2$ norm,*

$$\lim_{M\to\infty} \frac{F_2(\xi)F_4(\xi)}{F_3(\xi)}\chi(\xi)$$
$$= -\frac{1}{\sqrt{2\pi}}\exp\left\{-\frac{1}{4}\frac{\gamma^4 a^2 + \gamma^2 a^2 + 4\gamma^2 a + \gamma^2 + 1}{(\gamma^2 a + 1)^2}\eta^2 + \varepsilon\xi\right\}\chi(\xi).$$

The proof of Lemma 14 is obvious. The main ingredient in the proof of Lemma 15 is (64) and the fact that $F_4(\xi)/F_3(\xi)$ approaches to 1 uniformly on $[T,\infty)$.

We need another convergence result on $U(\xi)\widetilde{P}_{Ta}(\xi,\eta)U^{-1}(\eta)$:

PROPOSITION 8. *In trace norm,*

(70)
$$\lim_{M\to\infty} U(\xi)\widetilde{P}_{Ta}(\xi,\eta)U^{-1}(\eta) = 0.$$

The proof is left to the reader. Since all the four entries in $\widetilde{P}_{Ta}(\xi,\eta)$ can be expressed by Laguerre polynomials like (32) and (33), the asymptotic results like (63) and (64) give the convergence (70).

By Lemmas 14 and 15 and Proposition 8, we get in trace norm

$$\lim_{M\to\infty} \det(I - \widetilde{P}_T(\xi,\eta))$$
$$= \lim_{M\to\infty} \det(I - U(\xi)\widetilde{P}_T(\xi,\eta)U^{-1}(\eta))$$
$$= \lim_{M\to\infty} \det(I - U(\xi)\widetilde{P}_{Tb}(\xi,\eta)U^{-1}(\eta))$$
$$= \left(\int_\infty^T \frac{1}{\sqrt{2\pi}}e^{-t^2/2}\,dt\right)^2,$$

and we get the proof of the $a > \gamma^{-1}$ part of Theorem 1.

**4. Proof of $F_{\mathbf{GSE1}} = F_{\mathbf{GOE}}$.** In manipulation of kernels, we follow the method of [30]. The procedure seems informal and cursory, but is carefully justified in [30].

For notational simplicity, we denote $[\chi(\xi) = \chi_{(T,\infty)}(\xi)]$

$$B(\xi) = 1 - s^{(1)}(\xi) = \int_\xi^\infty \text{Ai}(t)\,dt.$$

First, we express the integral operator

$$\chi(\xi)\overline{\overline{P}}(\xi,\eta)\chi(\eta) = \begin{pmatrix} \chi(\xi)\overline{\overline{S}}_4(\xi,\eta)\chi(\eta) & \chi(\xi)\overline{\overline{SD}}_4(\xi,\eta)\chi(\eta) \\ \chi(\xi)\overline{\overline{IS}}_4(\xi,\eta)\chi(\eta) & \chi(\xi)\overline{\overline{S}}_4(\eta,\xi,)\chi(\eta) \end{pmatrix}$$

by
$$\begin{pmatrix} \chi(\xi)\dfrac{\partial}{\partial \xi} & 0 \\ 0 & \chi(\xi) \end{pmatrix} \begin{pmatrix} \overline{\overline{IS}}_4(\xi,\eta)\chi(\eta) & \overline{\overline{S}}_4(\eta,\xi)\chi(\eta) \\ \overline{\overline{IS}}_4(\xi,\eta)\chi(\eta) & \overline{\overline{S}}_4(\eta,\xi,)\chi(\eta) \end{pmatrix},$$

since by (22)–(23) and taking limit,
$$\frac{\partial}{\partial \xi}\overline{\overline{IS}}_4(\xi,\eta) = \overline{\overline{S}}_4(\xi,\eta),$$
$$\frac{\partial}{\partial \xi}\overline{\overline{S}}_4(\eta,\xi) = \overline{\overline{SD}}_4(\xi,\eta).$$

Then using (21) for $A$ bounded and $B$ trace class, upon suitably defining the Hilbert spaces our operators $A$ and $B$ are acting on, we find
$$\det\left(I - \begin{pmatrix} \chi(\xi)\dfrac{\partial}{\partial \xi} & 0 \\ 0 & \chi(\xi) \end{pmatrix}\begin{pmatrix} \overline{\overline{IS}}_4(\xi,\eta)\chi(\eta) & \overline{\overline{S}}_4(\eta,\xi)\chi(\eta) \\ \overline{\overline{IS}}_4(\xi,\eta)\chi(\eta) & \overline{\overline{S}}_4(\eta,\xi,)\chi(\eta) \end{pmatrix}\right)$$
$$= \det\left(I - \begin{pmatrix} \overline{\overline{IS}}_4(\xi,\eta)\chi(\eta) & \overline{\overline{S}}_4(\eta,\xi)\chi(\eta) \\ \overline{\overline{IS}}_4(\xi,\eta)\chi(\eta) & \overline{\overline{S}}_4(\eta,\xi,)\chi(\eta) \end{pmatrix}\begin{pmatrix} \chi(\eta)\dfrac{\partial}{\partial \eta} & 0 \\ 0 & \chi(\eta) \end{pmatrix}\right)$$
$$= \det\left(I - \begin{pmatrix} \overline{\overline{IS}}_4(\xi,\eta)\chi(\eta)\dfrac{\partial}{\partial \eta} & \overline{\overline{S}}_4(\eta,\xi)\chi(\eta) \\ \overline{\overline{IS}}_4(\xi,\eta)\chi(\eta)\dfrac{\partial}{\partial \eta} & \overline{\overline{S}}_4(\eta,\xi,)\chi(\eta) \end{pmatrix}\right),$$

and by conjugation with $\begin{pmatrix} 1 & 0 \\ -1 & 1 \end{pmatrix}$, we get
$$= \det\left(I - \begin{pmatrix} \overline{\overline{IS}}_4(\xi,\eta)\chi(\eta)\dfrac{\partial}{\partial \eta} + \overline{\overline{S}}_4(\eta,\xi)\chi(\eta) & \overline{\overline{S}}_4(\eta,\xi)\chi(\eta) \\ 0 & 0 \end{pmatrix}\right)$$
$$= \det\left(I - \left(\overline{\overline{IS}}_4(\xi,\eta)\chi(\eta)\dfrac{\partial}{\partial \eta} + \overline{\overline{S}}_4(\eta,\xi)\chi(\eta)\right)\right).$$

Since
$$\int_T^\infty \overline{\overline{IS}}_4(\xi,\eta)\frac{\partial}{\partial \eta}f(\eta)\,d\eta = \overline{\overline{IS}}_4(\xi,\eta)f(\eta)\big|_{\eta=T}^{\eta=\infty} - \int_T^\infty \frac{\partial}{\partial \eta}\overline{\overline{IS}}_4(\xi,\eta)f(\eta)\,d\eta,$$

as an operator
$$\overline{\overline{IS}}_4(\xi,\eta)\chi(\eta)\frac{\partial}{\partial \eta} = \overline{\overline{IS}}_4(\xi,\infty)\delta_\infty(\eta) - \overline{\overline{IS}}_4(\xi,T)\delta_T(\eta) - \frac{\partial}{\partial \eta}\overline{\overline{IS}}_4(\xi,\eta)\chi(\eta),$$

where $\delta_\infty$ and $\delta_T$ are (generalized) Dirac functions. Then with the help of identity
$$\int_\xi^\infty K_{\text{Airy}}(t,\eta)\,dt + \int_\eta^\infty K_{\text{Airy}}(\xi,t)\,dt = \int_\xi^\infty \text{Ai}(t)\,dt\int_\xi^\infty \text{Ai}(t)\,dt,$$

which can be proved directly from (1), we get

$$I - \left(\overline{\overline{IS}}_4(\xi,\eta)\chi(\eta)\frac{\partial}{\partial\eta} + \overline{\overline{S}}_4(\eta,\xi)\chi(\eta)\right)$$

$$= I - \left(K_{\text{Airy}}(\xi,\eta) - \frac{1}{2}B(\xi)\operatorname{Ai}(\eta) + \operatorname{Ai}(\eta)\right)\chi(\eta)$$

$$+ \left(\frac{1}{2}\int_T^\infty K_{\text{Airy}}(\xi,t)\,dt - \frac{1}{4}B(T)B(\xi) - \frac{1}{2}B(\xi) + \frac{1}{2}B(T)\right)\delta_T(\eta)$$

$$+ \frac{1}{2}B(\xi)\delta_\infty(\eta).$$

Now we denote $R(\xi,\eta)$ as the resolvent of $K_{\text{Airy}}(\xi,\eta)\chi(\eta)$, such that as integral operators

(71) $$I + R(\xi,\eta) = (I - K_{\text{Airy}}(\xi,\eta)\chi(\eta))^{-1},$$

then

$$I - \left(\overline{\overline{IS}}_4(\xi,\eta)\chi(\eta)\frac{\partial}{\partial\eta} + \overline{\overline{S}}_4(\eta,\xi)\chi(\eta)\right)$$

$$= (I - K_{\text{Airy}}(\xi,\eta)\chi(\eta))$$

$$\times \left(I - (I+R)\left(1 - \frac{1}{2}B(\xi)\right)\operatorname{Ai}(\eta)\chi(\eta)\right.$$

$$+ (I+R)\left(\frac{1}{2}\int_T^\infty K_{\text{Airy}}(\xi,t)\,dt\right.$$

$$\left. - \frac{1}{4}B(T)B(\xi) - \frac{1}{2}B(\xi) + \frac{1}{2}B(T)\right)\delta_T(\eta)$$

$$\left. + \frac{1}{2}(I+R)B(\xi)\delta_\infty(\eta)\right).$$

Again by the formula (21), in the form of (formula (17) in [30])

(72) $$\det\left(I - \sum_{k=1}^n \alpha_k \otimes \beta_k\right) = \det(\delta_{j,k} - (\alpha_j,\beta_k))_{j,k=1,\ldots,n}$$

we get

$$\det\left(I - (I+R)\left(1 - \frac{1}{2}B(\xi)\right)\operatorname{Ai}(\eta)\chi(\eta)\right.$$

$$+ (I+R)\left(\frac{1}{2}\int_T^\infty K_{\text{Airy}}(\xi,t)\,dt - \frac{1}{4}B(T)B(\xi) - \frac{1}{2}B(\xi) + \frac{1}{2}B(T)\right)\delta_T(\eta)$$

$$\left. + \frac{1}{2}(I+R)B(\xi)\delta_\infty(\eta)\right)$$

$$= \det \begin{pmatrix} 1+\alpha_{11} & \alpha_{12} & \alpha_{13} \\ \alpha_{21} & 1+\alpha_{22} & \alpha_{23} \\ \alpha_{31} & \alpha_{32} & 1+\alpha_{33} \end{pmatrix},$$

where upon the definition

$$\langle f(\xi), g(\xi) \rangle_T = \int_T^\infty f(\xi) g(\xi) \, d\xi,$$

we define

$$\alpha_{11} = \langle (I+R)(1 - \tfrac{1}{2}B(\xi)), -\operatorname{Ai}(\xi) \rangle_T,$$

$$\alpha_{12} = \Big\langle (I+R)\Big(\frac{1}{2}\int_T^\infty K_{\text{Airy}}(\xi, t)\, dt - \frac{1}{4}B(T)B(\xi) - \frac{1}{2}B(\xi) + \frac{1}{2}B(T)\Big),$$
$$- \operatorname{Ai}(\xi) \Big\rangle_T,$$

$$\alpha_{13} = \langle \tfrac{1}{2}(I+R)B(\xi), -\operatorname{Ai}(\xi)\rangle_T,$$

$$\alpha_{21} = (I+R)(1 - \tfrac{1}{2}B(\xi))|_{\xi=T},$$

$$\alpha_{22} = (I+R)\Big(\frac{1}{2}\int_T^\infty K_{\text{Airy}}(\xi, t)\, dt - \frac{1}{4}B(T)B(\xi) - \frac{1}{2}B(\xi) + \frac{1}{2}B(T)\Big)\Big|_{\xi=T},$$

$$\alpha_{23} = \tfrac{1}{2}(I+R)B(\xi)|_{\xi=T},$$

$$\alpha_{31} = (I+R)(1 - \tfrac{1}{2}B(\xi))|_{\xi=\infty} = 1,$$

$$\alpha_{32} = (I+R)\Big(\frac{1}{2}\int_T^\infty K_{\text{Airy}}(\xi, t)\, dt - \frac{1}{4}B(T)B(\xi) - \frac{1}{2}B(\xi) + \frac{1}{2}B(T)\Big)\Big|_{\xi=\infty}$$

$$= \frac{1}{2}B(T),$$

$$\alpha_{33} = \tfrac{1}{2}(I+R)B(\xi)|_{\xi=\infty} = 0.$$

If we take elementary row operations, we get

$$\det \begin{pmatrix} 1+\alpha_{11} & \alpha_{12} & \alpha_{13} \\ \alpha_{21} & 1+\alpha_{22} & \alpha_{23} \\ \alpha_{31} & \alpha_{32} & 1+\alpha_{33} \end{pmatrix}$$

$$= \det \begin{pmatrix} 1+\alpha_{11}-\alpha_{13} & \alpha_{12} - \tfrac{1}{2}B(T)\alpha_{13} & \alpha_{13} \\ \alpha_{21}-\alpha_{23} & 1+\alpha_{22} - \tfrac{1}{2}B(T)\alpha_{23} & \alpha_{23} \\ 0 & 0 & 1 \end{pmatrix}$$

$$= \det \begin{pmatrix} 1+\beta_{11} & \beta_{12} \\ \beta_{21} & 1+\beta_{22} \end{pmatrix},$$

where

$$\beta_{11} = \langle (I+R)(1-B(\xi)), -\operatorname{Ai}(\xi)\rangle_T,$$

$$\beta_{12} = \left\langle \frac{1}{2}(I+R)\left(\int_T^\infty K_{\text{Airy}}(\xi,t)\,dt - B(T)B(\xi) - B(\xi) + B(T)\right), -\text{Ai}(\xi) \right\rangle_T,$$

$$\beta_{21} = (I+R)(1-B(\xi))|_{\xi=T},$$

$$\beta_{22} = \frac{1}{2}(I+R)\left(\int_T^\infty K_{\text{Airy}}(\xi,t)\,dt - B(T)B(\xi) - B(\xi) + B(T)\right)\bigg|_{\xi=T}.$$

Using (71) and (72), we observe $[s^{(1)}(\xi) = 1 - B(\xi)]$

$$\det(I - K_{\text{Airy}}(\xi,\eta)\chi(\eta))\det\begin{pmatrix} 1+\beta_{11} & \beta_{12} \\ \beta_{21} & 1+\beta_{22} \end{pmatrix}$$

$$= \det\bigg(I - (K_{\text{Airy}}(\xi,\eta)\chi(\eta) + s^{(1)}(\xi)\text{Ai}(\eta))\chi(\eta)$$

$$+ \frac{1}{2}\bigg(\int_T^\infty K_{\text{Airy}}(\xi,t)\,dt - B(T)B(\xi) - B(\xi) + B(T)\bigg)\delta_T(\eta)\bigg).$$

If we denote $\tilde{R}(\xi,\eta)$ as the resolvent of $(K_{\text{Airy}}(\xi,\eta)\chi(\eta) + s^{(1)}(\xi)\text{Ai}(\eta))\chi(\eta)$, so that as operators

$$I + \tilde{R}(\xi,\eta) = (I + (K_{\text{Airy}}(\xi,\eta)\chi(\eta) + s^{(1)}(\xi)\text{Ai}(\eta))\chi(\eta))^{-1}$$

and

$$Q(\xi) = (I + \tilde{R})\bigg(\int_T^\infty K_{\text{Airy}}(\xi,t)\,dt - B(T)B(\xi) - B(\xi) + B(T)\bigg),$$

then

$$F_{\text{GSE1}} = \det(I - (K_{\text{Airy}}(\xi,\eta)\chi(\eta) + s^{(1)}(\xi)\text{Ai}(\eta))\chi(\eta))\det(I + \tfrac{1}{2}Q(\xi)\delta_T(\eta)).$$

To prove Theorem 2, we need only (6) and

$$\det(I + \tfrac{1}{2}Q(\xi)\delta_T(\eta)) = 1,$$

which by (72) is equivalent to

(73) $$Q(T) = 0.$$

If we take $f(\xi) = Q(\xi) + 1$, then (73) is

$$(I - (K_{\text{Airy}}(\xi,\eta)\chi(\eta) + s^{(1)}(\xi)\text{Ai}(\eta))\chi(\eta))(f(\xi) - 1)$$

$$= \int_T^\infty K_{\text{Airy}}(\xi,t)\,dt - B(T)B(\xi) - B(\xi) + B(T),$$

which is equivalent to

(74) $$(I - (K_{\text{Airy}}(\xi,\eta)\chi(\eta) + s^{(1)}(\xi)\text{Ai}(\eta))\chi(\eta))f(\xi) = s^{(1)}(\xi).$$

The integral equation (74) is solvable, and the solution is

$$f(\xi) = \frac{(I+R)s^{(1)}(\xi)}{1 - \langle (I+R)s^{(1)}(\xi), \text{Ai}(\xi)\rangle_T}.$$

Therefore to prove Theorem 2 we need only to prove $f(T) = 1$, which is equivalent to
$$(I+R)s^{(1)}(T) = 1 - \langle (I+R)s^{(1)}(\xi), \text{Ai}(\xi) \rangle_T.$$
This is a nontrivial result, but it can be derived by results in [30]. In Section VII of [30] Tracy and Widom define function $\bar{q}$ and $\bar{u}$ for both GOE and GSE. Our $(I+R)s^{(1)}(T)$ is equal to $\sqrt{2}$ times their $\bar{q}$ in GOE and our $\langle (I+R)s^{(1)}(\xi), \text{Ai}(\xi) \rangle_T$ is equal to 2 times their $\bar{u}$ in GOE. With

(75) $$(I+R)s^{(1)}(T) = e^{-\int_T^\infty q(s)\,ds},$$

(76) $$\langle (I+R)s^{(1)}(\xi), \text{Ai}(\xi) \rangle_T = 1 - e^{-\int_T^\infty q(s)\,ds},$$

where $q$ is the Painlevé II function determined by the differential equation
$$q''(s) = sq(s) + 2q^3(s)$$
together with the condition $q(s) \sim \text{Ai}(s)$ as $s \to \infty$.

We can give a proof of (75) and (76), based on the method and results in [29]. First, assume $T$ is fixed, then $(I+R)s^{(1)}$ is a function, and we have
$$\frac{d}{d\xi}(I+R)s^{(1)}(\xi) = (I+R)\frac{ds^{(1)}(\xi)}{d\xi} + \left[\frac{d}{d\xi}, (1+R)\right]s^{(1)}(\xi).$$
Since $\frac{d}{d\xi}s^{(1)}(\xi) = \text{Ai}(\xi)$ and we have (2.13) in [29], which is
$$\left[\frac{d}{d\xi}, (1+R)\right] = -(2+R)\text{Ai}(\xi) \cdot (1-K^t)^{-1}(\text{Ai}(\eta)\chi(\eta)) + R(\eta, T) \cdot \rho(T, \eta),$$
where $\rho(x,y) = \delta(x-y) + R(x,y)$ is the distribution kernel of $1+R$, and $K^t$ is the transpose (as an operator) of $K_{\text{Airy}}(\xi, \eta)\chi(\eta)$, we have
$$\frac{d}{d\xi}(I+R)s^{(1)}(\xi) = (1+R)\text{Ai}(\xi) - (1+R)\text{Ai}(\xi) \cdot \langle (I+R)s^{(1)}(\xi), \text{Ai}(\xi) \rangle_T$$
$$+ R(\xi, T) \cdot (1+R)s^{(1)}(T).$$
If we regard $T$ as a parameter, then we have

(77) $$\frac{d}{dT}(I+R)s^{(1)}(\xi; T) = -R(\xi, T) \cdot (1+R)s^{(1)}(T),$$

because (2.16) in [29] gives
$$\frac{1}{dT}(1+R) = R(\xi, T) \cdot \rho(T, \eta).$$
Therefore, if we set $\xi = T$ and take the derivative with respect to the parameter $T$, we have
$$\frac{d}{dT}((1+R)s^{(1)}(T)) = \left(\frac{d}{d\xi} + \frac{d}{dT}\right)((1+R)s^{(1)}(T))\bigg|_{\xi=T}$$
$$= (1+R)\text{Ai}(T) \cdot (1 - \langle (I+R)s^{(1)}(\xi), \text{Ai}(\xi) \rangle_T).$$

On the other hand, by (77) we have

$$\frac{d}{dT}\langle (I+R)s^{(1)}(\xi), \mathrm{Ai}(\xi)\rangle_T$$
$$= -(1+R)s^{(1)}(T) \cdot \mathrm{Ai}(T) + \left\langle \frac{d}{dT}(I+R)s^{(1)}(\xi), \mathrm{Ai}(\xi)\right\rangle_T$$
$$= -(1+R)s^{(1)}(T) \cdot \left(\mathrm{Ai}(T) + \int_T^\infty R(\xi,T)\,\mathrm{Ai}(\xi)\,d\xi\right)$$
$$= -(1+R)s^{(1)}(T) \cdot (1+R)\,\mathrm{Ai}(T).$$

(1.11) and (1.12) in [29] give the result

$$(1+R)\,\mathrm{Ai}(T) = q(T),$$

and now if we denote $(I+R)s^{(1)}(T) = s_T$ and $\langle (I+R)s^{(1)}(\xi), \mathrm{Ai}(\xi)\rangle_T = w_T$, we have

$$\begin{cases} \dfrac{d}{dT}s_T = q(1-w_T) \\ \dfrac{d}{dT}(1-w_T) = qs_T. \end{cases}$$

Now we can get (75) and (76) by boundary conditions.

## APPENDIX: DISCUSSION ON THE TRACE NORM CONVERGENCE OF INTEGRAL OPERATORS RELATED TO LUE

For convenience, we write (44) as

$$(78) \qquad L_j^{(2(M-N))}(2My) = \frac{e^{2My}}{2\pi i}\oint_C e^{-2Myz}\frac{z^{2(M-N)+j}}{(z-1)^{j+1}}\,dz,$$

where $C$ is a contour around 1, and we have another integral representation of Laguerre polynomials

$$(79) \qquad \begin{aligned} L_j^{(2(M-N))}(2Mx) &= \frac{(2(M-N)+j)!}{j!(2M)^{2(M-N)}}\frac{1}{x^{2(M-N)}2\pi i} \\ &\quad \times \oint_D e^{2Mxz}\frac{(z-1)^j}{z^{2(M-N)+j+1}}\,dz, \end{aligned}$$

where $D$ is a contour around 0.

Recall the integral operator $K(x,y)$ [11] for the rescaled LUE with parameters $2N$ and $2M$, and by (78) and (79) we have

$$K(x,y) = \sum_{j=0}^{2N-1} \frac{j!}{(2(M-N)+j)!}(2M)^{2(M-N)+1}$$

$$\times L_j^{(2(M-N))}(2Mx)L_j^{(2(M-N))}(2My)x^{M-N}y^{M-N}e^{-x+y}$$

(80)
$$= \frac{2M}{(2\pi i)^2} \frac{y^{M-N}e^{My}}{x^{M-N}e^{Mx}} \sum_{j=0}^{2N-1} \oint_C dz \oint_D dw \, e^{-2Myz} \frac{z^{2(M-N)+j}}{(z-1)^{j+1}}$$

$$\times e^{2Mxw} \frac{(w-1)^j}{w^{2(M-N)+j+1}}.$$

We can write the sum of integrands in (80) as

$$\sum_{j=0}^{2N-1} e^{2Mxz} \frac{z^{2(M-N)+j}}{(z-1)^{j+1}} e^{-2Myw} \frac{(w-1)^j}{w^{2(M-N)+j+1}}$$

$$= e^{2Mxz} e^{-2Myw} \frac{z^{2(M-N)}}{w^{2(M-N)}} \frac{1}{(z-1)w} \sum_{j=0}^{2N-1} \left(\frac{z(w-1)}{(z-1)w}\right)^j$$

$$= e^{2Mxz} e^{-2Myw} \frac{z^{2(M-N)}}{w^{2(M-N)}} \frac{1}{(z-1)w} \frac{1-(z(w-1)/((z-1)w))^{2N}}{1-z(w-1)/((z-1)w)}$$

$$= \frac{1}{z-w} e^{2Mxz} z^{2(M-N)} e^{-2Myw} \frac{1}{w^{2(M-N)}}$$

$$- \frac{1}{z-w} e^{2Mxz} \frac{z^{2M}}{(z-1)^{2N}} e^{-Myw} \frac{(w-1)^{2N}}{w^{2M}}.$$

By the residue theorem, let $C$ and $D$ be disjoint, then for the variable $z$, the pole $z = w$ is outside of $C$,

$$\oint_C dz \oint_D dw \frac{1}{z-w} e^{2Mxz} z^{2(M-N)} e^{-2Myw} \frac{1}{w^{2(M-N)}} = 0.$$

On the other side, we assume $\Re(w-z)$ to be less than 0, and get

$$\frac{1}{z-w} = 2M \int_0^\infty e^{t2M(w-z)} \, dt,$$

so that we have

$$\frac{2M}{(2\pi i)^2} \oint_C dz \oint_D dw \frac{1}{z-w} e^{-2Myz} \frac{z^{2M}}{(z-1)^{2N}} e^{2Mxw} \frac{(w-1)^{2N}}{w^{2M}}$$

(81)
$$= \frac{(2M)^2}{(2\pi i)^2} \oint_C dz \oint_D dw \int_0^\infty e^{-2M(y+t)z} \frac{z^{2M}}{(z-1)^{2N}}$$

$$\times e^{2M(x+t)w} \frac{(w-1)^{2N}}{w^{2M}}$$

$$= (2M)^2 \int_0^\infty \left(\frac{1}{2\pi i} \oint_C e^{-2M(y+t)z} \frac{z^{2M}}{(z-1)^{2N}} \, dz\right)$$

$$\times \left( \frac{1}{2\pi i} \oint_D e^{2M(x+t)w} \frac{(w-1)^{2N}}{w^{2M}} \, dw \right) dt.$$

Put (80)–(81) together, we get the result

$$K(x,y) = -(2M)^2 \frac{y^{M-N} e^{My}}{x^{M-N} e^{Mx}}$$

(82)
$$\times \int_0^\infty \left( \frac{1}{2\pi i} \oint_D e^{2M(x+t)w} \frac{(w-1)^{2N}}{w^{2M}} \, dw \right)$$

$$\times \left( \frac{1}{2\pi i} \oint_C e^{-2M(y+t)z} \frac{z^{2M}}{(z-1)^{2N}} \, dz \right) dt.$$

To find the probability that the largest eigenvalue $\geq T$ in the LUE, we need to consider the integral operator from $L^2([0,\infty))$ to $L^2([0,\infty))$ with the kernel $\chi(x)K(x,y)\chi(y)$. We can decompose it into the product of two integral operators by (82):

$$\chi(x)K(x,y)\chi(y) = -(2M)^2 \chi(x)J(x,t)\chi_{[0,\infty)}(t) \circ \chi_{[0,\infty)}(t)H(t,y)\chi(y),$$

where $\chi(x)J(x,t)\chi_{[0,\infty)}(t)$ and $\chi_{[0,\infty)}(t)H(t,y)\chi(y)$ stands for two integral operators with these kernels, and

$$J(x,t) = \frac{1}{x^{M-N} e^{Mx}} \frac{1}{2\pi i} \oint_D e^{2M(x+t)w} \frac{(w-1)^{2N}}{w^{2M}} \, dw,$$

$$H(t,y) = y^{M-N} e^{My} \frac{1}{2\pi i} \oint_C e^{-2M(y+t)z} \frac{z^{2M}}{(z-1)^{2N}} \, dz.$$

Since we consider the limiting distribution of the largest eigenvalue around $(1+\gamma^{-1})^2$, we take $p = (1+\gamma^{-1})^2$, $q = \frac{(1+\gamma)^{4/3}}{\gamma(2M)^{2/3}}$, $x = p + q\xi$, $y = p + q\eta$ and $t = q\tau$. Then for the rescaled kernel $\chi(\xi)\widetilde{K}(\xi,\eta)\chi(\eta)$, we have

$$\chi(\xi)\widetilde{K}(\xi,\eta)\chi(\eta) = \chi(\xi)\widetilde{J}(\xi,\tau)\chi_{[0,\infty)}(\tau) \circ \chi_{[0,\infty)}(\tau)\widetilde{H}(\tau,\eta)\chi(\eta),$$

where

$$\widetilde{J}(\xi,\tau) = \frac{(\gamma+1)^{4/3}}{\gamma} M^{1/3} \frac{\gamma^{2N} e^{N-M}}{x^{M-N} e^{Mx}} \frac{1}{2\pi i} \oint_D e^{2M(p+q(\xi+\tau))w} \frac{(w-1)^{2N}}{w^{2M}} \, dw,$$

$$\widetilde{H}(\tau,\eta) = \frac{(\gamma+1)^{4/3}}{\gamma} M^{1/3} \frac{y^{M-N} e^{My}}{\gamma^{2N} e^{N-M}} \frac{1}{2\pi i} \oint_C e^{-2M(p+q(\xi+\tau))z} \frac{z^{2M}}{(z-1)^{2N}} \, dz.$$

We want to prove the trace norm convergence

(83)
$$\lim_{M \to \infty} \chi(\xi)\widetilde{K}(\xi,\eta)\chi(\eta)$$
$$= \chi(\xi)K_{\text{Airy}}(\xi,\eta)\chi(\eta)$$
$$= \chi(\xi)\operatorname{Ai}(\xi+\tau)\chi_{[0,\infty)}(\tau) \circ \chi_{[0,\infty)}(\tau)\operatorname{Ai}(\tau+\eta)\chi(\eta).$$

By results in functional analysis, we need only to prove the convergence in Hilbert–Schmidt norm of (e.g., see [20])

$$\lim_{M \to \infty} \chi(\xi)\widetilde{J}(\xi,\tau)\chi_{[0,\infty)}(\tau) = \chi(\xi)K_{f\text{Airy}}(\xi,\eta)\chi(\eta), \tag{84}$$

$$\lim_{M \to \infty} \chi_{[0,\infty)}(\tau)\widetilde{H}(\tau,\eta)\chi(\eta) = \chi_{[0,\infty)}(\tau)\text{Ai}(\tau+\eta)\chi(\eta). \tag{85}$$

Since for integral operators, the convergence in Hilbert–Schmidt norm is equivalent to the convergence in $L^2$ norm of their kernels as two variable functions, we can verify (84) and (85) by asymptotic analysis similar to that in Section 3.

For the integral operator $\chi(\xi)\widetilde{S}_{4a1}(\xi,\eta)\chi(\eta)$ in (39), we have

$$\chi(\xi)\widetilde{S}_{4a1}(\xi,\eta)\chi(\eta) = \frac{1}{2}\chi(\xi)\frac{1}{\sqrt{x}}\widetilde{\widetilde{J}}(\xi,\tau)\chi_{[0,\infty)}(\tau) \circ \chi_{[0,\infty)}(\tau)\sqrt{y}\widetilde{\widetilde{H}}(\tau,\eta)\chi(\eta),$$

where we define $\widetilde{\widetilde{J}}$ and $\widetilde{\widetilde{H}}$ in the same way as $\widetilde{J}$ and $\widetilde{J}$, but use parameters $2N-2$ and $2M-2$ instead of $2N$ and $2M$. Similarly, in the $\widetilde{S}_{4a1}$ part of (61), we have

$$\chi(\xi)e^{\varepsilon\xi}\widetilde{S}_{4a1}(\xi,\eta)e^{-\varepsilon\eta}\chi(\eta)$$
$$= \frac{1}{2}\chi(\xi)\frac{e^{\varepsilon\xi}}{\sqrt{x}}\widetilde{\widetilde{J}}(\xi,\tau)\chi_{[0,\infty)}(\tau) \circ \chi_{[0,\infty)}(\tau)\frac{\sqrt{y}}{e^{-\varepsilon\eta}}\widetilde{\widetilde{H}}(\tau,\eta)\chi(\eta).$$

We can give rigorous proofs to (39) and (61) in the same way as (83).

**Acknowledgments.** The author is most grateful to his advisor Mark Adler, who pointed out the spiked model problem and gave warm encouragement to the author. I also thank Ira Gessel for help in combinatorics, especially his suggestions in proofs in Section 2, and Jinho Baik for his suggestions on related models.

DEPARTMENT OF MATHEMATICS
BRANDEIS UNIVERSITY
WALTHAM, MASSACHUSETTS 02454
USA
E-MAIL: wangdong@brandeis.edu